\input amstex
\magnification=1200
\documentstyle{gen-j}
\hsize=15truecm
\vsize=23truecm
\hfuzz 4pt
\document

\def\lft#1{{#1}\hfill}
\def\hquad{\hskip.5em{}}
\def\vl{\vrule height 9.5pt depth 6.5pt}
\font\sevenrm=cmr7
\def\dsp{\displaystyle}
\def\hl{\leaders\hrule height 2.8pt depth -2.5pt\hfill}
\def\lra{\hbox{\raise 0.3335pt\hbox{$\longrightarrow$}}}
\def\llra{\hbox to 30pt{\hl\hskip -3pt$\lra$}}
\def\lllra{\hbox to 40pt{\hl\hskip -3pt$\lra$}}
\def\llllra{\hbox to 80pt{\hl\hskip -3pt$\lra$}}
\def\ora{\mathop{\lllra}}
\def\orasm{\mathop{\llra}}
\def\lratwo{\mathop{\llllra}}
\def\lram{\mathop{\longrightarrow}}

\font\sf=cmss10
\font\titl=cmr10 scaled\magstep2
\font\tit=cmr10 scaled\magstep3
\font\ti=cmr9 scaled\magstep2

\font\bigse=cmbsy10 scaled\magstep 5
\def\largese{\hbox{\bigse\char'045}}  
 5
%
\font\bigne=cmbsy10 scaled\magstep 5
\def\largene{\hbox{\bigne\char'045}}  
\font\bigse=cmbsy10 scaled\magstep 5
\def\largese{\hbox{\bigse\char'046}}  
\catcode`\@=11
\def\bqed{\ifhmode\unskip\nobreak\fi\quad
  \ifmmode\blacksquare\else$\m@th\blacksquare$\fi}
\def\mon{\qopname@{monotonely}}
\catcode`\@=\active
\nologo
\topmatter
\endtopmatter
\centerline{{\tit Shuffle Invariance of the}}\vskip0.2cm
\centerline{{\tit Super-RSK Algorithm}}
\vskip0.8cm
\centerline{{\sf Amitai\ Regev\footnote{Department of Theoretical
Mathematics, Weizmann Institute of Science, Rehovot
76100, Israel;\quad {\tt regev,tamars\@wisdom.weizmann.ac.il}\hfil\break
{\tt http://www.wisdom.weizmann.ac.il/\~{}regev,tamars}},\quad\quad Tamar
Seeman${}^1$}}
\vskip0.9cm
{\bf Abstract} As in the $(k,l)$-RSK (Robinson-Schensted-Knuth)
of [1], other super-RSK algorithms can be
applied to sequences of variables from the
set \hfill\break
$\{t_1,...,t_k,u_1,...,u_l\}$, where $t_1<\cdots<t_k$, and
$u_1<\cdots<u_l$. While the $(k,l)$-RSK of [1] is the case where
$t_i<u_j$ for all $i$ and $j$, these other super-RSK's correspond to
all the
$\big{(}{{k+l}\atop{k}}\big{)}$
shuffles of the $t$'s and $u$'s satisfying the above
restrictions that $t_1<\cdots<t_k$ and $u_1<\cdots<u_l$.
We show that the shape of the tableaux produced by any such super-RSK is
independent of the particular shuffle of the $t$'s and $u$'s.
\bigskip

\noindent{\titl 1\ \ Introduction}\bigskip
We follow the tableaux-terminology of [7].
The classical Frobenius-Schur-Weyl theory shows how the
SSYT (Semi-Standard-Young-Tableaux) determine the representations of
$GL(m,\Bbb C)$ (or $gl(m,\Bbb C)$). Here $GL(m,\Bbb C)$ ($gl(m,\Bbb C)$)
is the General Linear Lie group (algebra).
Also, SYT (Standard-Young-Tableaux) play an important role here.
The notion of $(k,l)$ SSYT 
is introduced  in [1], where similar relationships between such tableaux
and the representations of $pl(k,l)$ are shown. Here  $pl(k,l)$
is the General Linear Lie super-algebra.

The $(k,l)$ SSYT are defined, via a $(k,l)$-RSK algorithm, as follows [1].
Fix integers $k,\/ l\ge 0$,
$k+l>0$, and $k+l$ symbols
$t_1,...,t_k,u_1,...,u_l$ such that $t_1<...<t_k<u_1<...<u_l$. Let
$$a_{k,l}(n)=\Big{\{}\Big{(}{{1...n}\atop{v_1...v_n}}\Big{)}\Big{|}\ v_i\in
\{t_1,...,t_k,u_1,...u_l\}\Big{\}}.$$

\noindent To map $a_{k,l}(n)$ to pairs of tableaux $(P,Q)$, apply to
each $v\in a_{k,l}(n)$ the $(k,l)$-RSK, in which the usual
RSK insertion algorithm [7] is applied to the $t_i$'s, and the
conjugate correspondence (see [1]) is applied to the $u_j$'s; see
the examples below. By the definitions of [1],
the {\it insertion
tableau}, $P=P(v)$, mapped from $v\in a_{k,l}(n),$ is $(k,l)$ semistandard;
that
is, it satisfies the following three properties:
\roster
\item "(a)" The ``$t$ part'' (i.e., the cells filled with $t_i$'s) is a
tableau.
\item "(b)" The $t_i$'s are nondecreasing in rows, strictly increasing in
columns.
\item "(c)" The $u_j$'s are nondecreasing in columns, strictly increasing in
rows.
\endroster\medskip

As in the usual correspondence, the {\it recording tableau}, $Q=Q(v)$,
indicates the
order in which the new cells were added to $P$. Clearly, $Q$ is SYT having 
the same shape as that of $P.$
\medskip
A total order of $\{t_1,...,t_k,u_1,...,u_l\}$ which is compatible with
$t_1<\cdots <t_k$ and $u_1<\cdots <u_l$, is called a {\it shuffle} (of
$t_1,...,t_k$ and $u_1,...,u_l$). For example,
$t_1 < u_1 < u_2 < t_2$ is such a shuffle, compatible with
$t_1 < t_2$ and $u_1 < u_2.$
Clearly, there are
$\big{(}{{k+l}\atop{k}}\big{)}$ such shuffles; of these, Berele and
Regev chose to work with
$t_1<\cdots <t_k<u_1<\cdots <u_l$, which we call {\it the} $(k,l)$
{\it shuffle} (see [1, 2.4]).
The shuffle $t_1 < u_1 < t_2 < u_2 < \cdots < t_k < u_k,$ with its
corresponding SSYT, appears in section 4 of [3].

Let $I=I(k,l)$ denote the set of all such $\big{(}{{k+l}\atop{k}}\big{)}$
shuffles. Given $A\in I$, there is a corresponding $A$-RSK insertion
algorithm; if
$v\in a_{k,l}(n)$, then $v{\dsp\lram_A}(P,Q)$ by that algorithm.
$P=P_A=P(v,A)$ is the insertion tableau, and
$Q=Q_A=Q(v,A)$ is the recording tableau. Here $P$ is an $A$-SSYT;
that is, it satisfies the following three properties.\medskip
\roster
\item "(a)" $P$ is weakly $A$-increasing in both rows and columns.
\item "(b)" The $t_i$'s are strictly increasing in columns.
\item "(c)" The $u_j$'s are strictly increasing in rows.
\endroster\bigskip

\demo{\bf Example} Let $k=l=2$, $A,B \in I=I(2,2),$ where

$$A: \ t_1 < t_2 < u_1 < u_2 \quad\text{and}\quad B: \ u_1 < u_2 < t_1 <
t_2.$$
Let $$v=\Big{(}{{1\cdots\cdots 4}\atop{u_2,t_1,t_2,u_1}}\Big{)}.$$
\medskip
Then
\overfullrule=0pt
$$ v\lram_A\  \vcenter{\tabskip=0pt\offinterlineskip\halign
{\lft{$#$}&&\hquad\lft{$#$}\hquad&\lft{$#$}\cr
\multispan3{\hrulefill}\cr
\vl& u_2 &\vl \cr
\multispan3{\hrulefill}\cr}}\quad\quad
\vcenter{\tabskip=0pt\offinterlineskip\halign
{\lft{$#$}&&\hquad\lft{$#$}\hquad&\lft{$#$}\cr
\multispan5{\hrulefill}\cr
\vl& t_1 &\vl& u_2 &\vl \cr
\multispan5{\hrulefill}\cr}}\quad\quad
\vcenter{\tabskip=0pt\offinterlineskip\halign
{\lft{$#$}&&\hquad\lft{$#$}\hquad&\lft{$#$}\cr
\multispan7{\hrulefill}\cr
\vl& t_1 &\vl& t_2 &\vl& u_2 &\vl \cr
\multispan7{\hrulefill}\cr}}\quad\quad
\vcenter{\tabskip=0pt\offinterlineskip\halign
{\lft{$#$}&&\hquad\lft{$#$}\hquad&\lft{$#$}\cr
\multispan7{\hrulefill}\cr
\vl&t_1&\vl&t_2&\vl&u_2&\vl \cr
\multispan7{\hrulefill}\cr
\vl&u_1&\vl \cr
\multispan3{\hrulefill}\cr}}\  =\ P_A,\quad\quad\text{while}$$

\overfullrule=0pt
$$ v\lram_B\  \vcenter{\tabskip=0pt\offinterlineskip\halign
{\lft{$#$}&&\hquad\lft{$#$}\hquad&\lft{$#$}\cr
\multispan3{\hrulefill}\cr
\vl&u_2&\vl \cr
\multispan3{\hrulefill}\cr}}\quad\quad
\vcenter{\tabskip=0pt\offinterlineskip\halign
{\lft{$#$}&&\hquad\lft{$#$}\hquad&\lft{$#$}\cr
\multispan5{\hrulefill}\cr
\vl&u_2&\vl&t_1&\vl \cr
\multispan5{\hrulefill}\cr}}\quad\quad
\vcenter{\tabskip=0pt\offinterlineskip\halign
{\lft{$#$}&&\hquad\lft{$#$}\hquad&\lft{$#$}\cr
\multispan7{\hrulefill}\cr
\vl&u_2&\vl&t_1&\vl&t_2&\vl \cr
\multispan7{\hrulefill}\cr}}\quad\quad
\vcenter{\tabskip=0pt\offinterlineskip\halign
{\lft{$#$}&&\hquad\lft{$#$}\hquad&\lft{$#$}\cr
\multispan7{\hrulefill}\cr
\vl&u_1&\vl&u_2&\vl&t_2&\vl \cr
\multispan7{\hrulefill}\cr
\vl&t_1&\vl \cr
\multispan3{\hrulefill}\cr}}\  =\ P_B$$\medskip

Thus $v{\dsp\lram_A}(P_A,Q)\ $ and $\ v{\dsp\lram_B}(P_B,Q),$
where
\overfullrule=0pt
$$ Q\ =\  \vcenter{\tabskip=0pt\offinterlineskip\halign
{\lft{$#$}&&\hquad\lft{$#$}\hquad&\lft{$#$}\cr
\multispan7{\hrulefill}\cr
\vl&1&\vl&2&\vl&3&\vl \cr
\multispan7{\hrulefill}\cr
\vl&4&\vl \cr
\multispan3{\hrulefill}\cr}}, \quad\text{and }P_A\text{ and }P_B\text{ are
as above.}\bqed$$
\enddemo\bigskip

\definition{\bf Definition} Denote by sh$(v,A)=$ sh$(P_A)$ the shape of
the insertion tableau $P(v,A)=P_A$ of $v\in a_{k,l}(n)$ under the $A$-RSK.
\enddefinition\bigskip

Given a shuffle $A \in I$ and the pair $(P,Q),$ where $P$ is
$A$-SSYT, $Q$ is  SYT, and $sh(P)=sh(Q)$, the $A$ insertion
algorithm can obviously be reversed. By standard arguments (see
for example [7, chap. 7]) this yields \medskip

\proclaim{\bf Theorem~1} Let $A \in I$ be a shuffle.
Then the $A$-RSK insertion
algorithm $v{\dsp\lram_A}(P_A,Q_A)$ is a bijection between
$a_{k,l}(n)$ and
$$
\{(P_A,Q_A) \ |\ P_A\ \text{is}\ A\text{\rm -SSYT,} \ Q_A \ \text{is } 
\text{\rm SYT,\ \ sh}(P_A)\ =\ \text{\rm sh}(Q_A) \}.
$$\endproclaim\bigskip

\definition{\bf Remark}
Denote such a tableau $P=(P_{i,j})$ and denote $<_A$ by $<$. Clearly, if

\noindent  $P_{i,j}=t_r$ then $P_{i,j-1}\le
P_{i,j}\le P_{i,j+1}$ and $P_{i-1,j}<P{i,j}<P_{i,j+1}$. Similarly, 

\noindent if $P_{i,j}=u_r$ then $P_{i,j-1}<P_{i,j}<P_{i,j+1}$ and
$P_{i-1,j}\le P{i,j}\le P_{i,j+1}$.
\enddefinition\bigskip

Denote by sh($v,A$) the shape of tableaux $P(v,A)$ and $Q(v,A)$.
This brings us to our main result.\bigskip

\proclaim{\bf Theorem~2}
Let $v\in a_{k,l}(n)$,  $A,B\in I$, $v{\dsp\lram_A}(P_A,Q_A)$
and $v{\dsp\lram_B}(P_B,Q_B)$. Then {\rm sh}$(P_A)=\ ${\rm sh}$(P_B)$. 
Consequently, $Q_A=Q_B$.
\endproclaim\bigskip

In other words, the shape of the tableau obtained through
any of the $(k,l)$-shuffle-RSK algorithms,
is independent of the particular shuffle of the $t$'s and $u$'s.\bigskip

\definition{\bf Definition}
Let $A\in I$ and $\lambda\vdash n$,
i.e. a partition of $n$. Let $\Im_A(\lambda)$ denote the set of
the
$A$-SSYT of shape $\lambda$:
$$\Im_A(\lambda)=\{T\ |\ T\ \text{is } A\text{-SSYT, sh}(T)=\lambda\}.$$
\enddefinition\bigskip

Recall the definition of $type(T)$ from [7, page 309].

Theorem 2 implies 
\medskip

\proclaim{\bf Theorem 3 [6]}
Let $A,B\in I$, $\lambda\vdash n$. Then there exists a bijection $\varphi:
\Im_A(\lambda)\rightarrow \Im_B(\lambda)$ such that for all $\ T\in
\Im_A(\lambda)$,
type$(T)=$ type$(\varphi(T))$. (In fact, there exist (at least) $d_\lambda$
such canonical bijections, where $d_\lambda$ is the number of SYT's of shape
$\lambda$.)
\endproclaim\medskip

Theorem~3 appears in [6], where it is proven by a different method. Our proof
of the theorem is as follows.\bigskip

\demo{\bf Proof of Theorem~3}
Is based on the following diagram:

$$\aligned &\quad\quad\quad (P_A,Q)\\
&\magnification=1200
\largene\hbox{\hskip -12pt{$_{ A\text{-RSK}}$}}\\
v\in a_{k,l}(n),\quad v& \\
&\largese\hbox{\hskip -12pt{$^{ B\text{-RSK}}$}}\\
&\quad\quad\quad (P_B,Q)\endaligned$$\medskip

Thus choose a SYT $Q$ of shape $\lambda$. Given $P=P_A\in \Im_A(\lambda)$,
get
\def\tog{\hbox{\rm inverse\lower 12pt\hbox{\hskip -30pt 
\rm A-RSK}}}
$$ (P_A,Q)~\ora_{\tog}~~~~v~~~~
\ora_{\hbox{B-RSK}}~~~~(P_B,Q)~. $$
\medskip

\noindent This defines the bijection $\varphi=\varphi_{_Q}:\ \varphi(P_A)=P_B$. 
Clearly, type$(P_A)=$ type$(P_B)$ and by Theorem 2,
sh$(P_A)=$ sh$(P_B).$
\bqed\enddemo\bigskip

Recall from [2] the notation $w(T)$ for the weight of a tableau $T$. For
example, let \medskip
\overfullrule=0pt
$$ T\ =\ \vcenter{\tabskip=0pt\offinterlineskip\halign
{\lft{$#$}&&\hquad\lft{$#$}\hquad&\lft{$#$}\cr
\multispan9{\hrulefill}\cr
\vl&t_1&\vl&t_1&\vl&u_2&\vl&u_3&\vl \cr
\multispan9{\hrulefill}\cr
\vl&t_2&\vl&t_3&\vl&u_2&\vl \cr
\multispan7{\hrulefill}\cr
\vl&u_1&\vl&u_3&\vl \cr
\multispan5{\hrulefill}\cr
\vl&u_1&\vl \cr
\multispan3{\hrulefill}\cr}}$$
\medskip

\noindent then $w(T)=x_1^2x_2x_3y_1^2y_2^2y_3^2$. Also recall the ``hook''
(or the ``super'') Schur function

$$HS_\lambda(x;y)=HS_\lambda(x_1,...,x_k;y_1,...,y_l)\quad\text{[1],[2]}.$$
\medskip

When $A$ is the shuffle $A_0:\ t_1<\cdots <t_k<u_1<\cdots 
<u_l$, $HS_{\lambda}(x;y)$ is given by
$$ HS_\lambda(x_1,...,x_k;y_1,...,y_l)=\sum_{T\in \Im_{A_0}
(\lambda)}w(T) $$
[1, Thm. 6.10]. See also [4], [5] and [6].\medskip
It clearly  follows from Theorem~3 that\medskip

\proclaim{\bf Corollary~4}
For any $A\in I$, 
$$ HS_\lambda(x_1,...,x_k;y_1,...,y_l)=
\sum_{T\in \Im_A(\lambda)}w(T)~.$$

\endproclaim\bigskip

Given a shuffle $A\in I$, the $A$-RSK is based on $A$, on the regular RSK 
for the $t_i$'s and the conjugate-regular RSK for the $u_j$'s. 

In addition to the regular RSK, there is also the dual RSK [7, page 331]. Given
the shuffle $A\in I$, this leads to four possible $A$-insertion algorithms: 
either the regular or the dual for the $t_i$'s, and either the conjugate
regular or the conjugate dual for the $u_j$'s. In fact, the previous $A$-RSK
is: ($t$-regular, $u$-conjugate-regular), which we denote as the 
(regular, regular)-$A$-RSK. Similarly, ($t$-regular, $u$-dual-conjugate) is the
(regular, dual)-$A$-RSK. Similarly for the algorithms (dual, regular)-$A$-RSK 
and (dual, dual)-$A$-RSK. Each of these three new insertion algorithms
exhibits a similar shape invariance under all shuffles $A\in I$.\bigskip

\proclaim{\bf Theorem~5}
\roster 
\item "(a)" Let $v\in a_{k,l}(n)$, and $A,B\in I$ such that 
$$ v~\lratwo_{\hbox{\rm(regular,regular)-A-RSK}}~(P^*_A,Q^*_A),\quad\quad v~
\lratwo_{\hbox{\rm(regular,regular)-B-RSK}}~(P^*_B,Q^*_B).$$ 
Then {\rm sh}$(P^*_A)=$
{\rm sh}$(P^*_B)$. Consequently, $Q^*_A=Q^*_B$.\medskip

\item "(b)" Let $v\in a_{k,l}(n)$, and $A,B\in I$ such that 
$$ v~\lratwo_{\hbox{\rm(regular,dual)-A-RSK}}~(P^*_A,Q^*_A),\quad\quad v~
\lratwo_{\hbox{\rm(regular,dual)-B-RSK}}~(P^*_B,Q^*_B).$$ 
Then {\rm sh}$(P^*_A)=$
{\rm sh}$(P^*_B)$. Consequently, $Q^*_A=Q^*_B$.\medskip
\item "(c)" Let $v\in a_{k,l}(n)$, and $A,B\in I$ such that 
$$ v~\lratwo_{\hbox{\rm(dual,regular)-A-RSK}}~(P^*_A,Q^*_A),\quad\quad v~ 
\lratwo_{\hbox{\rm(dual,regular)-B-RSK}}~(P^*_B,Q^*_B).$$ 
Then {\rm sh}$(P^*_A)=$
{\rm sh}$(P^*_B)$. Consequently, $Q^*_A=Q^*_B$.\medskip
\item "(d)" Let $v\in a_{k,l}(n)$, and $A,B\in I$ such that 
$$v~\lratwo_{\hbox{\rm(dual,dual)-A-RSK}}~(P^*_A,Q^*_A),\quad\quad v~ 
\lratwo_{\hbox{\rm(dual,dual)-B-RSK}}~(P^*_B,Q^*_B).$$ 
Then {\rm sh}$(P^*_A)=$
{\rm sh}$(P^*_B)$. Consequently, $Q^*_A=Q^*_B$.
\endroster\endproclaim\bigskip

Clearly, Theorem 5(a) is Theorem 2 above.
The proof of Theorem~2 is given in the next section, 
which is the main body of this paper. First we describe
the $A$-RSK algorithm in details. The main step in the proof of Theorem~2 is
Lemma~2.15. It shows that a transposition of the variables in the shuffle
(i.e a single change in the order of some $t_i$ and $u_j$), does not alter
the shape of the resulting tableaux. 
In section 3 we prove the remaining parts (b), (c), and (d) of Theorem 5,
essentially by deducing them from Theorem 2.
\bigskip

\noindent{\titl 2\ \ Invariance of Shape}\bigskip

As in the $(k,l)$-RSK, the $A$-RSK insertion algorithm involves applying
the usual RSK correspondence to the $t_i$'s, and the
conjugate correspondence to the $u_j$'s. This is illustrated in the
following example.\bigskip

\definition{\bf Definition~2.1}
For $i,j\in\Bbb Z^+$, let $c(i,j)$ denote the cell in row $i$ and column $j$
of a given tableau.
\enddefinition\bigskip

\demo{\bf Example~2.2}
Under the shuffle $A=t_1<u_1<t_2<u_2<t_3$, perform the insertion
\overfullrule=0pt
$$ \vcenter{\tabskip=0pt\offinterlineskip\halign
{\lft{$#$}&&\hquad\lft{$#$}\hquad&\lft{$#$}\cr
\multispan7{\hrulefill}\cr
\vl&u_1&\vl&t_2&\vl&t_2&\vl \cr
\multispan7{\hrulefill}\cr
\vl&u_1&\vl&u_2&\vl \cr
\multispan5{\hrulefill}\cr
\vl&t_3&\vl \cr
\multispan3{\hrulefill}\cr}}~\leftarrow t_1. $$
\bigskip

(a) $t_1<u_1$ $\implies$ $t_1$ occupies $c(1,1)$. Now, a $u_i$ is always
bumped to the next column, hence $u_1$ is bumped to column $2$.

(b) $u_1<t_2$ $\implies$ $u_1$ occupies $c(1,2)$. Now, a $t_i$ is always
bumped to the next row, hence $t_2$ is bumped to row $2$.

(c) $u_1<t_2<u_2$ $\implies$ $t_2$ occupies $c(2,2)$, bumping $u_2$ to
column
$3$.

(d) $u_2>t_2$ $\implies$ $u_2$ settles in $c(2,3)$.

\overfullrule=0pt
$$ \text{(a)}\quad \vcenter{\tabskip=0pt\offinterlineskip\halign
{\lft{$#$}&&\hquad\lft{$#$}\hquad&\lft{$#$}\cr
\multispan7{\hrulefill}\cr
\vl&t_1&\vl&t_2&\vl&t_2&\vl \cr
\multispan7{\hrulefill}\cr
\vl&u_1&\vl&u_2&\vl \cr
\multispan5{\hrulefill}\cr
\vl&t_3&\vl \cr
\multispan3{\hrulefill}\cr}}\quad,\quad\quad
\text{(b)}\quad \vcenter{\tabskip=0pt\offinterlineskip\halign
{\lft{$#$}&&\hquad\lft{$#$}\hquad&\lft{$#$}\cr
\multispan7{\hrulefill}\cr
\vl&t_1&\vl&u_1&\vl&t_2&\vl \cr
\multispan7{\hrulefill}\cr
\vl&u_1&\vl&u_2&\vl \cr
\multispan5{\hrulefill}\cr
\vl&t_3&\vl \cr
\multispan3{\hrulefill}\cr}}$$

$$ \text{(c)}\quad \vcenter{\tabskip=0pt\offinterlineskip\halign
{\lft{$#$}&&\hquad\lft{$#$}\hquad&\lft{$#$}\cr
\multispan7{\hrulefill}\cr
\vl&t_1&\vl&t_2&\vl&t_2&\vl \cr
\multispan7{\hrulefill}\cr
\vl&u_1&\vl&t_2&\vl \cr
\multispan5{\hrulefill}\cr
\vl&t_3&\vl \cr
\multispan3{\hrulefill}\cr}}\quad,\quad\quad
\text{(d)}\quad \vcenter{\tabskip=0pt\offinterlineskip\halign
{\lft{$#$}&&\hquad\lft{$#$}\hquad&\lft{$#$}\cr
\multispan7{\hrulefill}\cr
\vl&t_1&\vl&u_1&\vl&t_2&\vl \cr
\multispan7{\hrulefill}\cr
\vl&u_1&\vl&t_2&\vl&u_2&\vl \cr
\multispan7{\hrulefill}\cr
\vl&t_3&\vl \cr
\multispan3{\hrulefill}\cr}} $$
\bqed\enddemo\bigskip

The proof of Theorem~2 will follow from the following analysis of the
$A$-RSK
algorithm.\bigskip

\proclaim{\bf Lemma~2.3}
Let $p$ be an $A$-SSYT, $v\in\{t_1,...,t_k,u_1,...,u_l\}$. The insertion
$P\leftarrow v$ is made of a sequence of several steps. In an intermediate
$m$-th such a step, we have an $A$-SSYT $\widetilde{P}$ together with an
element $P_{i,j}$ that was bumped from $c(i,j)$ by ${\widetilde{P}}_{i,j}$,
${\widetilde{P}}_{i,j}{<\atop A}P_{i,j}$, and we need to do the following
insertion:\roster
\item "(a)" If $P_{i,j}=t_r$, insert it into the $i+1$-th row of
$\widetilde{P}$.
\item "(b)" If $P_{i,j}=u_s$, insert it into the $j+1$-th column of
$\widetilde{P}$.
\endroster
We show that in both cases, the result would be an $A$-SSYT $P^*$, and
-- except for the last step -- together with a new element
${\widetilde{P}}_{i',j'}$ (bumped from $c(i',j')$), which is to
be inserted into $P^*$. Moreover,\roster
\item "(1)" If $P_{i,j}=t_r$ then $c(i',j')=c(i+1,j')$ and $j'\le j$.
\item "(2)" If $P_{i,j}=u_s$ then $c(i',j')=c(i',j+1)$ and $i'\le i$.
\endroster
then
\endproclaim\bigskip

\demo{\bf Proof}
Note that (2) is obtained from (1) by conjugation, hence it suffices to just
prove (1).

{\it Proof of} (1): Denote the $i$-th row of $\widetilde{P}$ by
$$a_1\cdots\cdots\cdots a_{j-1}{\widetilde{P}}_{i,j}a_{j+1}
\cdots\cdots\cdots a_g,$$ so
$a_j=P_{i,j}$ and by assumption, $P_{i,j}=t_r$. Thus
$$\aligned &\ \vdots\\ &a_1\cdots\cdots\cdots\cdots\cdots\cdots
a_{j-1}{\widetilde{P}}_{i,j}
a_{j+1}\cdots\cdots\cdots a_g\\ \widetilde{P}\quad=\quad
&b_1\cdots\cdots\cdots\cdots\cdots\cdots\cdots\cdots\cdots\cdots\cdots\cdots
\cdots b_f\\
&c_1\cdots\cdots\cdots\cdots\cdots\cdots\cdots\cdots\cdots\cdots\cdots\cdots
c_h\\ &\
\vdots\endaligned$$

\noindent
and $P_{i,j}=t_r$ is inserted into the $i+1$-th row $b_1\cdots\cdots b_f$.

Let $b_{j'-1}\le P_{i,j}<b_{i,j'}$, so in $P^*$, the $i+1$-th row is
$$b_1\cdots\cdots\cdots\cdots
b_{j'-1}P_{i,j}b_{j'+1}\cdots\cdots\cdots\cdots b_f.$$
Since ${\widetilde{P}}_{i,j}$ bumped $P_{i,j}$, we have
${\widetilde{P}}_{i,j}
<P_{i,j}$. Since $a_j=P_{i,j}=t_r$, hence $P_{i,j}<b_j$. Together with
$b_{j'-1}
\le P_{i,j}<b_{j'}$, this implies that $j'\le j$, hence
$$\aligned &\ \vdots\\ &a_1\cdots\cdots a_{j'-1}a_{j'}a_{j'+1}\cdots\cdots
{\widetilde{P}}_{ij}a_{j+1}\cdots\cdots\cdots a_g\\P^*\quad=\quad
&b_1\cdots\cdots b_{j'-1}P_{ij}b_{j'+1}\cdots\cdots\cdot
b_jb_{j+1}\cdots\cdots
\cdots
b_f\\ &c_1\cdots\cdots c_{j'-1}c_{j'}c_{j'+1}\cdots\cdots\cdot\cdot
c_jc_{j+1}\cdots\cdots c_h\\ &\ \vdots\endaligned$$

By the induction assumption on $\widetilde{P}$, we only need to verify that
the part
$$\aligned &a_{j'}\\ &P_{i,j}\\ &c_{j'}\endaligned$$\medskip
\noindent
of the $j'$-th column is $A$-semistandard, i.e.:
since $P_{i,j}=t_r$, we need to show that $a_{j'}\le{\widetilde{P}}_{i,j}
<c_{j'}$.
This follows from $a_{j'}\le{\widetilde{P}}_{i,j}<P_{i,j}=t_r<b_{j'}
\le c_{j'}$.
\bqed\enddemo\bigskip

\definition{\bf Definition~2.4}
Two shuffles $A,B\in I$ are {\it adjacent} if there exist $t_i$
and $u_j$ such that\roster
\item "1)" $t_i<u_j$ in $A$.
\item "2)" $u_j<t_i$ in $B$.
\item "3)" All other pairs have the same order relations in $A$ and in $B$.
\endroster
\noindent In that case, call $A$ and $B$ $(t_i,u_j)$-{\it adjacent}. Thus
$A$ and $B$ differ by the transposition $(t_i,u_j)$.
\enddefinition\bigskip

\definition{\bf Remark~2.5}
Trivially, for any $A$,$B\in I$ there exist $A_0$, $A_1$,...,$A_n\in I$ such
that $A_0=A$, $A_n=B$, and $A_r$ is adjacent to $A_{r+1}$, $0\le r\le n-1$.
Thus to prove Theorem~1, it suffices to show that for all
$v\in a_{k,l}(n)$ and for every pair $(A,B)$ of adjacent shuffles,
sh$(v,A)=$ sh$(v,B)$. Therefore for the rest of this section, let $A,B\in I$
be $(t_i,u_j)$-adjacent, with $t_i<_{_A}u_j$ and $u_j<_{_B}t_i$.
\enddefinition\bigskip

\proclaim{\bf Lemma~2.6}
Let $A\in I$, $w\in a_{k,l}(n)$, and for some
$x\in\{t_1,...,t_k,u_1,...,u_l\}$, let $w'$ be the sequence obtained by
omitting from $w$ all elements $A$-greater than $x$.
Let $P_A$ and $P_A'$ be the insertion
tableaux obtained from $w$ and $w'$ respectively under
shuffle $A$. Then $P_A'$ is a subtableau of $P_A$.
\endproclaim\bigskip

\demo{\bf Proof}
Let $w~{\dsp\orasm_{\hbox{\sevenrm A-RSK}}}~P_A;\ P:\emptyset, 
P_1, P_2,...,P_n=P_A$, and similarly let $w'~
{\dsp\orasm_{\hbox{\sevenrm A-RSK}}}~P'_A$; $P':
\emptyset, P'_1, P'_2,...,P'_m=P'_A\quad (m=|w'|)$.

Assume $P'_i$ is a subtableau of $P_{j_i}$, and insert (a corresponding) $y$
in $w$.

If $x<_{_A}y$, $y$ is not in $w'$ so $P'_i$ is not affected. Also,
inserting $y$ into $P_{j_i}$, $y$ does not affect the subtableau $P'_i
\subseteq P_{j_i}$, since $y$ bumps only elements that are $A$-greater than
itself.

A similar argument applies when $y\le x$: now $y$ is also in $w'$, and is
inserted into $P'_i$ and into $P_{j_i}$. Clearly, in $P_{j_i}$ it is also
inserted into the subtableau $P'_i\subseteq P_{j_i}$, and the proof follows.
\bqed\enddemo\bigskip

\proclaim{\bf Corollary~2.7}
Let $A,B\in I$ be $(t_i,u_j)$-adjacent, $v\in a_{k,l}(n)$, $v{\dsp\lram
_A}(P_A,Q_A)$ and $v{\dsp\lram_B}(P_B,Q_B)$. Then the elements that
are
both $A$-less and $B$-less than $t_i$ and $u_j$ form identical subtableaux
in $P_A$ and $P_B$.
\endproclaim\bigskip

\demo{\bf Proof}
Denote by $v'$ the sequence obtained by omitting from $v$ all elements
($A$- and $B$-) greater
than or equal to $t_i$ and $u_j$. By $(t_i,u_j)$-adjacency, the largest
element smaller than $t_i$ and $u_j$, in both $A$ and $B$, is the same
element
$x$. Moreover, $v'$ is obtained by omitting from $v$ all elements which are
($A$- or $B$-) greater than $x$. Let $P_A'$ and $P_B'$ denote the insertion
tableaux of $v'$ under shuffles $A$ and $B$ respectively. Then by Lemma~2.6,
$P_A'$ and $P_B'$ are subtableaux of $P_A$ and $P_B$ respectively. But the
elements that are $A$- or $B$-less than $t_i$ and $u_j$ are ordered
identically in $A$ and $B$, so $P_A'=P_B'$.
\bqed\enddemo\bigskip

\definition{\bf Notation}
As above, let $A,B\in I$ be two shuffles that are $(t_i,u_j)$-adjacent:
$t_i<u_j$ in $A$ and $u_j<t_i$ in $B$. Let
$v\in a_{k,l}(n)$, and denote $v{\dsp\lram_A}(P_A,Q_A)$ and
$v{\dsp\lram_B}(P_B,Q_B)$.
\enddefinition\bigskip

\definition{\bf Notation}
Given the tableau $P_A$ (and similarly for $P_B$), let regions 1, 2 and 3
denote, respectively, the
regions occupied (1) by elements less than $t_i$ and $u_j$, (2) by $t_i$ and
$u_j$, and (3) by elements greater than $t_i$ and $u_j$.
\enddefinition\bigskip

\demo{\bf Example~2.8}
Let $v=u_1t_3t_2u_2t_2u_1t_1$, and let
$$A=t_1<u_1<t_2<u_2<t_3\atop B=t_1<u_1<u_2<t_2<t_3.$$
Then $A$ and $B$ are $(t_i,u_j)$-adjacent, with $t_i=t_2$ and $u_j=u_2$, and
\medskip
\overfullrule=0pt
$$P_A\ =\ \vcenter{\tabskip=0pt\offinterlineskip\halign
{\lft{$#$}&&\hquad\lft{$#$}\hquad&\lft{$#$}\cr
\multispan7{\hrulefill}\cr
\vl&t_1&\vl&u_1&\vl&t_2&\vl \cr
\multispan7{\hrulefill}\cr
\vl&u_1&\vl&t_2&\vl&u_2&\vl \cr
\multispan7{\hrulefill}\cr
\vl&t_3&\vl \cr
\multispan3{\hrulefill}\cr}}\ ,\quad\quad
P_B\ =\ \vcenter{\tabskip=0pt\offinterlineskip\halign
{\lft{$#$}&&\hquad\lft{$#$}\hquad&\lft{$#$}\cr
\multispan7{\hrulefill}\cr
\vl&t_1&\vl&u_1&\vl&u_2&\vl \cr
\multispan7{\hrulefill}\cr
\vl&u_1&\vl&t_2&\vl&t_2&\vl \cr
\multispan7{\hrulefill}\cr
\vl&t_3&\vl \cr
\multispan3{\hrulefill}\cr}}\ .$$\bigskip
In both tableaux, region 1 contains the elements $t_1$ and $u_1$, region 2
contains $t_2$ and $u_2$, and region 3 contains $t_3$. Note that in this
example, regions 1
and 3 are the same in $P_A$ as in $P_B$, and region 2 is identically shaped
in $P_A$ and $P_B$. We shall show that this is always true.
\bqed\enddemo\bigskip

By Lemma~2.6, both region 1, as well as the union of regions 1 and 2, form
subtableaux in $P$. It is easy to check that region 2 does not contain the
configuration
\overfullrule=0pt
$$ \vcenter{\tabskip=0pt\offinterlineskip\halign
{\lft{$#$}&&\hquad\lft{$#$}\hquad&\lft{$#$}\cr
\multispan5{\hrulefill}\cr
\vl&a&\vl&b&\vl \cr
\multispan5{\hrulefill}\cr
\vl&c&\vl&d&\vl \cr
\multispan5{\hrulefill}\cr}}~. $$

If it does, assume $d=t_i$. Then $b=u_j$, so
$u_j<t_i$, and $a\ne t_i,u_j$. Similarly if $d=u_j$. It follows that region
2
forms part of the rim of the subtableaux which is the union of regions 1 and
2.\bigskip

\definition{\bf Remark~2.9}
Note that (part of) region 2 in $P_A$ (i.e. $t_i<u_j$) always looks like

$${\aligned&\quad\quad\quad\quad\ t_i\cdots\cdots t_i\\
&\quad\quad\quad\quad\ u_j\\
&\quad\quad\quad\quad\ \vdots\\ &t_i\cdots\cdots t_iu_j\\
&u_j\\ &\vdots\\ &u_j\endaligned}$$

Namely: Except possibly for the rightmost element, all other elements in a
row
are $t_i$'s. Similarly, except for possibly the top element, all other
elements in a column are $u_j$'s.

Similarly, in $P_B$ (i.e. $u_j<t_i$), part of region 2 looks like

$${\aligned&\quad\quad\quad\quad\ u_jt_i\cdots\cdots t_i\\
&\quad\quad\quad\quad\ \vdots\\ &u_jt_i\cdots\cdots t_i\\
&\vdots\\ &u_j\endaligned}$$
\enddefinition\bigskip

Denote $v=v_1\cdots v_n$. The tableau $P_A$ is created by applying the
$A$-RSK
insertion algorithm to each of $v_1,...,v_n$ successively. For each $v_m$,
let
$l_{m(A)}$ denote the length of the insertion path [7, page~317] of $v_m$
under shuffle $A$ -- that is, the number of insertion steps that occur when
$v_m$ is inserted while forming $P_A$.
The total number of insertion steps involved in the formation of $P_A$ is
thus
$s_A=\sum_{m=1}^n l_{m(A)}$. For every $r\in\{1,...,s_A\}$, let $P_A^r$
be the insertion tableau as it appears immediately after insertion step $r$.

Similarly, under shuffle $B$, the length of the insertion path of $v_m$ into
$P_B$ is $l_{m(B)}$, and the total number of insertion steps involved in
forming
$P_B$ is $s_B=\sum_{m=1}^n l_{m(B)}$, with $P_B^r$ denoting the insertion
tableau after insertion step $r$.\bigskip

\demo{\bf Example~2.10}
As in Example~2.8, let $v=v_1\cdots v_7=u_1t_3t_2u_2t_2u_1t_1$, and let
$A=t_1<u_1<t_2<u_2<t_3$. Then tableau $P_A$ is formed by the $A$-RSK as
follows (ignore the underlines).\medskip
\overfullrule=0pt
$$ \vcenter{\tabskip=0pt\offinterlineskip\halign
{\lft{$#$}&&\hquad\lft{$#$}\hquad&\lft{$#$}\cr
\multispan3{\hrulefill}\cr
\vl&\underline{u_1}&\vl \cr
\multispan3{\hrulefill}\cr}}\quad\quad
\vcenter{\tabskip=0pt\offinterlineskip\halign
{\lft{$#$}&&\hquad\lft{$#$}\hquad&\lft{$#$}\cr
\multispan5{\hrulefill}\cr
\vl&u_1&\vl&\underline{t_3}&\vl \cr
\multispan5{\hrulefill}\cr}}\quad\quad
\vcenter{\tabskip=0pt\offinterlineskip\halign
{\lft{$#$}&&\hquad\lft{$#$}\hquad&\lft{$#$}\cr
\multispan5{\hrulefill}\cr
\vl&u_1&\vl&\underline{t_2}&\vl \cr
\multispan5{\hrulefill}\cr
\vl&\underline{t_3}&\vl \cr
\multispan3{\hrulefill}\cr}}\quad\quad
\vcenter{\tabskip=0pt\offinterlineskip\halign
{\lft{$#$}&&\hquad\lft{$#$}\hquad&\lft{$#$}\cr
\multispan5{\hrulefill}\cr
\vl&u_1&\vl&t_2&\vl \cr
\multispan5{\hrulefill}\cr
\vl&\underline{u_2}&\vl \cr
\multispan3{\hrulefill}\cr
\vl&\underline{t_3}&\vl \cr
\multispan3{\hrulefill}\cr}} $$

\overfullrule=0pt
$$ \vcenter{\tabskip=0pt\offinterlineskip\halign
{\lft{$#$}&&\hquad\lft{$#$}\hquad&\lft{$#$}\cr
\multispan7{\hrulefill}\cr
\vl&u_1&\vl&t_2&\vl&\underline{t_2}&\vl \cr
\multispan7{\hrulefill}\cr
\vl&u_2&\vl \cr
\multispan3{\hrulefill}\cr
\vl&t_3&\vl \cr
\multispan3{\hrulefill}\cr}} \quad\quad
\vcenter{\tabskip=0pt\offinterlineskip\halign
{\lft{$#$}&&\hquad\lft{$#$}\hquad&\lft{$#$}\cr
\multispan7{\hrulefill}\cr
\vl&u_1&\vl&t_2&\vl&t_2&\vl \cr
\multispan7{\hrulefill}\cr
\vl&\underline{u_1}&\vl&\underline{u_2}&\vl \cr
\multispan5{\hrulefill}\cr
\vl&t_3&\vl \cr
\multispan3{\hrulefill}\cr}} \quad\quad
\vcenter{\tabskip=0pt\offinterlineskip\halign
{\lft{$#$}&&\hquad\lft{$#$}\hquad&\lft{$#$}\cr
\multispan7{\hrulefill}\cr
\vl&\underline{t_1}&\vl&\underline{u_1}&\vl&t_2&\vl \cr
\multispan7{\hrulefill}\cr
\vl&u_1&\vl&\underline{t_2}&\vl&\underline{u_2}&\vl \cr
\multispan7{\hrulefill}\cr
\vl&t_3&\vl \cr
\multispan3{\hrulefill}\cr}} $$\medskip

For all $i\in\{1,...,7\}$, the underlined elements in tableau $i$ lie in the
insertion path of element $v_i$. Thus $l_{1(A)}=l_{2(A)}=l_{5(A)}=1$,
$l_{3(A)}
=l_{4(A)}=l_{6(A)}=2$, $l_{7(A)}=4$, and $s_A=\sum_{i=1}^7l_{i(A)}=13$.
If for example, $r=7=\sum_{i=1}^5l_{i(A)}$, then we have
\overfullrule=0pt
$$P^r\ =\ \vcenter{\tabskip=0pt\offinterlineskip\halign
{\lft{$#$}&&\hquad\lft{$#$}\hquad&\lft{$#$}\cr
\multispan7{\hrulefill}\cr
\vl&u_1&\vl&t_2&\vl&t_2&\vl \cr
\multispan7{\hrulefill}\cr
\vl&u_2&\vl \cr
\multispan3{\hrulefill}\cr
\vl&t_3&\vl \cr
\multispan3{\hrulefill}\cr}}\quad,\quad\quad
P^{r+1}\ =\ \vcenter{\tabskip=0pt\offinterlineskip\halign
{\lft{$#$}&&\hquad\lft{$#$}\hquad&\lft{$#$}\cr
\multispan7{\hrulefill}\cr
\vl&u_1&\vl&t_2&\vl&t_2&\vl \cr
\multispan7{\hrulefill}\cr
\vl&u_1&\vl \cr
\multispan3{\hrulefill}\cr
\vl&t_3&\vl \cr
\multispan3{\hrulefill}\cr}}\quad\quad\bqed$$
\enddemo\bigskip

\demo{\bf Example~2.11}
Let $k=l=1,\ \ A: t<u,\ \ B: u<t,\ \ v=v_1v_2=tu$. Then
\overfullrule=0pt
$$P_A\ :\ \emptyset,\quad \vcenter{\tabskip=0pt\offinterlineskip\halign
{\lft{$#$}&&\hquad\lft{$#$}\hquad&\lft{$#$}\cr
\multispan3{\hrulefill}\cr
\vl&\underline{t}&\vl \cr
\multispan3{\hrulefill}\cr}}\ ,\quad
\vcenter{\tabskip=0pt\offinterlineskip\halign
{\lft{$#$}&&\hquad\lft{$#$}\hquad&\lft{$#$}\cr
\multispan3{\hrulefill}\cr
\vl&t&\vl \cr
\multispan3{\hrulefill}\cr
\vl&\underline{u}&\vl \cr
\multispan3{\hrulefill}\cr}}\ ;\quad
l_{1(A)}=l_{2(A)}=1,\quad\quad\quad\ $$
$$P_B\ :\ \emptyset,\quad \vcenter{\tabskip=0pt\offinterlineskip\halign
{\lft{$#$}&&\hquad\lft{$#$}\hquad&\lft{$#$}\cr
\multispan3{\hrulefill}\cr
\vl&\underline{t}&\vl \cr
\multispan3{\hrulefill}\cr}}\ ,\quad
\vcenter{\tabskip=0pt\offinterlineskip\halign
{\lft{$#$}&&\hquad\lft{$#$}\hquad&\lft{$#$}\cr
\multispan3{\hrulefill}\cr
\vl&\underline{u}&\vl \cr
\multispan3{\hrulefill}\cr
\vl&\underline{t}&\vl \cr
\multispan3{\hrulefill}\cr}}\ ;\quad
l_{1(B)}=1,\ l_{2(B)}=2.\ \bqed$$
\enddemo\bigskip

\definition{\bf Definition~2.12}
For $p,q\in\Bbb Z^+$, we say that $P_A^p\sim P_B^q$ (with respect to the
formations of $P_A$ and $P_B$) if:
\roster
\item Regions 1 and 3 are identical in $P_A^p$ and $P_B^q$.
\item Region 2 is identically shaped in $P_A^p$ and $P_B^q$; moreover, in
each
connected component of that region 2, the number of $t_i$'s (hence of
$u_j$'s)
in $P_A^p$ equals the number of $t_i$'s (hence of $u_j$'s) in $P_B^q$.
\item Either $p=s_A$ and $q=s_B$, or both $p<s_A$ and $q<s_B$. In the latter
case, the next insertion step involves
inserting the same element into the same row (or column) in both tableaux.
\endroster\enddefinition\bigskip

\demo{\bf Example~2.13}
The tableaux of Example~2.8 satisfy $P_A\sim P_B$. Regions 1 and 3 in the
two
tableaux are identical, satisfying property 1 of Definition~2.12. Region 2
consists of one component which is identically shaped, and contains exactly
one $t_i$ and one $u_j$, in both tableaux. This verifies property 2. Since
both tableaux correspond to $p=s_A$ and $q=s_B$, property 3 is satisfied as
well.
\bqed\enddemo\bigskip

\proclaim{\bf Lemma~2.14}
For any shuffle $A\in I$, and for all $p\in\{2,...,s_A\}$ and
$r,s\in\Bbb Z^+$, if $c(r,s)$ contains some $w$ in $P_A^{p-1}$, then
$c(r,s)$
contains some $z\le_{_A}w$ in $P_A^p$.

Conversely, if $c(r,s)$ contains some element $z$ in $P_A^p$, then $c(r,s)$
was either empty or contained some $w\ge_{_A}z$ in $P_A^{p-1}$.
\endproclaim\bigskip

\demo{\bf Proof}
Follows from the $A$-RSK algorithm.
\bqed\enddemo\bigskip

\noindent{\bf The Proof of Theorem~2} clearly follows from\bigskip

\proclaim{\bf Lemma~2.15}
Let $A,B\in I$ be $(t_i,u_j)$-adjacent, $v\in a_{k,l}(n)$, 
$v~{\dsp\orasm_{\hbox{\sevenrm A-RSK}}}~(P_A,Q_A)$ and
$v~{\dsp\orasm_{\hbox{\sevenrm B-RSK}}}~(P_B,Q_B)$,
then $P_A\sim P_B$.
\endproclaim

\demo{\bf Proof}
We prove that $P_A\sim P_B$, by induction on the insertion steps of $P_A$
and
$P_B$. Trivially $P_A^1=P_B^1$. Now let $p\in\{1,...,s_A-1\}$,
$q\in\{1,...,s_B-1\}$ and assume that 1) $P_A^p\sim P_B^q$, and also 2)
either
$P_A^{p-1}\sim P_B^{q-1}$ or $P_A^{p-1}\sim P_B^{q-2}$ or $P_A^{p-2}\sim
P_B^{q-1}$. We show that this implies that either $P_A^{p+1}\sim P_B^{q+1}$
or $P_A^{p+1}\sim P_B^{q+2}$ or $P_A^{p+2}\sim P_B^{q+1}$. This clearly
implies
the proof of the lemma (by induction on $p+q$).

Note that if $P_A^p\sim P_B^q$, then by 2.12.3, step $p+1$ in $P_A$ and step
$q+1$ in
$P_B$ are identical; that is, the same element, $x$, is inserted into the
same row (or column) in both tableaux. We
assume that that $x$ is a $t$-element, and therefore enters some row,
denoted {\it
row} $r$; the case where $x$ is a $u$-element is analogous. Since $P_A^p\sim
P_B^q$, row $r$
is empty in $P_A^p$ if and only if it is empty in $P_B^q$. The case where
row $r$ is empty is trivial, so we assume throughout that row $r$ is
nonempty
in $P_A^p$ and $P_B^q$.
\medskip
\roster
\item Suppose that under both shuffles $A$ and $B$, $x>t_i$ and $u_j$. Since
$P_A^p\sim P_B^q$, the last nonempty
cell in row $r$ must be in the same region in both $P_A^p$ and $P_B^q$, and
if
it is in region 3, then it must be occupied by the same element in both
tableaux.\medskip

Case 1.1: Row $r$ in $P_A^p$ (and in $P_B^q$) terminates with an element
less
than or equal to $x$. In this case, $x$ is affixed to the end of the row in
both tableaux, so $P_A^{p+1}$ and $P_B^{q+1}$
have the same shape and clearly satisfy properties 1 and 2 of
Definition~2.12. Let $m$
denote the size of $P_A^{p+1}$ and $P_B^{q+1}$.
If $m=n$, which is the size of $P_A$ and $P_B$, then the insertion algorithm
terminates here. Otherwise, the next step is to begin $v_{m+1}$'s insertion
path by inserting $v_{m+1}$ into either the first row or the first column in
both tableaux. This verifies 2.12.3 and we have $P_A^{p+1}\sim P_B^{q+1}$.
\medskip

Case 1.2: Row $r$ in $P_A^p$ contains an element $z>x$ (under both $A$ and
$B$). Since $P_A^p\sim P_B^q$, the same is true in $P_B^q$. In this case,
$x$ bumps an element greater than itself --
a region-3 element -- and occupies its cell in both tableaux. Thus both the
cell occupied by $x$ and the element bumped by $x$ are
identical in the two tableaux, which 2.12.3. Since 2.12.1 and 2.12.2
clearly hold, it follows that $P_A^{p+1}\sim P_B^{q+1}$.\medskip
\item
Suppose that $x=t_i$. During step $P_B^q\rightarrow P_B^{q+1}$,
$x=t_i>_{_B}u_j$ bumps the first region-3 element in row
$r$, or if no such element exists, $x$ occupies the first empty cell in
that row. Let $c(r,s)$ be the cell occupied by $x$ in $P_B^{q+1}$.\medskip

Case 2.1: In row $r$ of $P_A^p$, region 2 either terminates with $t_i$ or
does
not
appear at all in that row. Then $x$ occupies $c(r,s)$ also in $P_A^{p+1}$
(and
bumps the same element as in $P_B^{q+1}$), so $P_A^{p+1}\sim P_B^{q+1}$.
\medskip

Case 2.2: In $P_A^p$, the last region-2 element in row $r$ is
$u_j$. Let this $u_j$ be in $c(r,s')$. Since $P_A^p\sim P_B^q$,
$c(r,s')$ is the last region-2 cell in row $r$ in both tableaux. Since in
$P_B^q\rightarrow P_B^{q+1}$, $x$ was inserted into $c(r,s)$,
we have $s=s'+1$.
Thus $u_j$ is in $c(r,s-1)$ and is bumped by $x=t_i$ to
column $s$ during $P_A^p\rightarrow P_A^{p+1}$.
We prove that in such a case, $P_A^{p+2}\sim P_B^{q+1}$. To do so, we
show that

2.2.1: In $P_A^{p+1}\rightarrow P_A^{p+2}$, $u_j$ settles in $c(r,s)$, to
the
immediate right of $x$.

2.2.2: This implies that~2.12.2 for $P_A^{p+2}\sim P_B^{q+1}$ is satisfied.

2.2.3: Both~2.12.1 and~2.12.3 for $P_A^{p+2}\sim P_B^{q+1}$ are satisfied.
\medskip

{\it Proof of 2.2.1:} If $r=1$, then $u_j$ clearly settles in $c(r,s)$ in
$P_A^{p+2}$. We therefore assume that $r>1$.\medskip
To prove that $u_j$ settles in $c(r,s)$ in $P_A^{p+2}$, we need only to show
that $c(r-1,s)$ in $P_A^{p+1}$ contains an element $b\le u_j$, since
$c(r,s)$
in $P_A^{p+1}$ contains some element $z>_{_A}u_j$.
Now, since $r>1$, $x=t_i$ arrived at row $r$ in $P_A^p$ (and similarly in
$P_B^q$)
after being bumped from row $r-1$ of $P_A^{p-1}$. Let $c(r-1,h)$ be the cell
occupied by $x$ in $P_A^{p-1}$, before it was bumped from row $r-1$.\medskip
\def\togone{\hbox{\sevenrm $x$ is bumped\lower 12pt\hbox{\hskip -45pt 
\sevenrm from $c(r-1,h)$}}}
\def\togtwo{\hbox{\sevenrm $x$ is inserted\lower 12pt\hbox{\hskip -45pt 
\sevenrm into $c(r,s-1)$}}}
\def\togthr{\hbox{\sevenrm $u_j$ is inserted\lower 12pt\hbox{\hskip -45pt 
\sevenrm into column $s$}}}
$$ P_A^{p-1}~\ora_{\togone}~P_A^p~\ora_{\togtwo}~
P_A^{p+1}~\ora_{\togthr}~P_A^{p+2} $$\medskip

Since $x$ is inserted into $c(r,s-1)$ of $P_A^{p+1}$,
Lemma~2.3 implies that $s-1\le h$. If $s\le h$,
then $c(r-1,s)$ was occupied by an element less than or equal to $x$ in
$P_A^{p-1}$, and this continues to be true in $P_A^p$ and $P_A^{p+1}$, which
implies our claim (that $b\le u_j$). On the
other hand, suppose that $h=s-1$. In this case, we prove that $c(r-1,s)$
contains an element less than or equal to $u_j$ by showing that otherwise
we would have a contradiction to the induction assumptions. Our proof of
this
is illustrated by the following figures, each of which consists of the
block of cells in rows $r-1$ and $r$ and columns $s-1$ and $s$ in the
corresponding tableaux.
\overfullrule=0pt
$$\quad P_A^{p-1}: \vcenter{\tabskip=0pt\offinterlineskip\halign
{\lft{$#$}&&\hquad\lft{$#$}\hquad&\lft{$#$}\cr
&&&&&$s$ \cr
&&&&&$\vdots$ \cr
&&\multispan5{\hrulefill}\cr
&&\vl&x&\vl&b&\vl \cr
&&\multispan5{\hrulefill}\cr
r&\cdots&\vl&&\vl&&\vl&, \cr
&&\multispan5{\hrulefill}\cr}}~~
P_A^p:\ \vcenter{\tabskip=0pt\offinterlineskip\halign
{\lft{$#$}&&\hquad\lft{$#$}\hquad&\lft{$#$}\cr
&&&&& s \cr
&&&&& \vdots \cr
&&\multispan5{\hrulefill}\cr
&&\vl&&\vl& b &\vl \cr
&&\multispan5{\hrulefill}\cr
r & \cdots &\vl& u_j &\vl&&\vl&, \cr
&&\multispan5{\hrulefill}\cr}} $$
$$ P_A^{p+1}:\ \vcenter{\tabskip=0pt\offinterlineskip\halign
{\lft{$#$}&&\hquad\lft{$#$}\hquad&\lft{$#$}\cr
&&&&& s \cr
&&&&& \vdots \cr
&&\multispan5{\hrulefill}\cr
&&\vl&&\vl& b &\vl \cr
&&\multispan5{\hrulefill}\cr
r &\cdots&\vl& x&\vl&&\vl \cr
&&\multispan5{\hrulefill}\cr}}~ $$
$$ P_B^{q-2}:\ \vcenter{\tabskip=0pt\offinterlineskip\halign
{\lft{$#$}&&\hquad\lft{$#$}\hquad&\lft{$#$}\cr
&&&&& s \cr
&&&&& \vdots \cr
&&\multispan5{\hrulefill}\cr
&&\vl&\ &\vl& z &\vl \cr
&&\multispan5{\hrulefill}\cr
r & \cdots &\vl&\ &\vl&&\vl&, \cr
&&\multispan5{\hrulefill}\cr}}~
P_B^{q-1}: \vcenter{\tabskip=0pt\offinterlineskip\halign
{\lft{$#$}&&\hquad\lft{$#$}\hquad&\lft{$#$}\cr
&&&&&$s$ \cr
&&&&& \vdots \cr
&&\multispan5{\hrulefill}\cr
&&\vl&\ &\vl& e &\vl \cr
&&\multispan5{\hrulefill}\cr
r & \cdots &\vl&&\vl&&\vl \cr
&&\multispan5{\hrulefill}\cr}}~. $$

So, assume that in $P_A^{p-1}$, $x$ occupied $c(r-1,s-1)$, and $c(r-1,s)$
was
either empty or contained an element $b>u_j$ (a region-3 element).
We show that this contradicts the claim of the induction hypothesis, that
either $P_A^{p-1}\sim P_B^{q-1}$, $P_A^{p-1}\sim P_B^{q-2}$, or
$P_A^{p-2}\sim P_B^{q-1}$.\medskip
\def\togone{\hbox{\sevenrm $e$ is inserted\lower 12pt\hbox{\hskip -45pt 
\sevenrm into $c(r-1,s)$}}}
\def\togtwo{\hbox{\sevenrm $x$ is bumped\lower 12pt\hbox{\hskip -45pt 
\sevenrm from $c(r-1,s')$}}}
\def\togthr{\hbox{\sevenrm $x$ is inserted\lower 12pt\hbox{\hskip -45pt 
\sevenrm into $c(r,s)$}}}
$$ P_B^{q-2}~\ora_{\togone}~P_B^{q-1}~\ora_{\togtwo}~
P_B^q~\ora_{\togthr}~
P_B^{q+1} $$\medskip
Now in $P_B^{q+1}$, $x=t_i$ occupies $c(r,s)$,
so $c(r-1,s)$ is occupied by an element less than $x$. Lemma~2.3 implies
that
in $P_B^q$, $x$ occupied $c(r-1,s')$, where $s\le s'$, and this implies that
in $P_B^{q-1}$, $c(r-1,s)$ contained an element $e\le x$, a region-1 or -2
element. But the same cell in $P_A^{p-1}$ was either empty or contained a
region-3 element, and by Lemma~2.14, the same must have been true in
$P_A^{p-2}$. Thus the above assumption, that $h=s-1$, implies that neither
$P_A^{p-1}\sim P_B^{q-1}$ nor $P_A^{p-2}\sim P_B^{q-1}$ is satisfied. We
show
now that it also implies that $P_A^{p-1}\nsim P_B^{q-2}$.

Suppose that $P_A^{p-1}\sim P_B^{q-2}$ is satisfied. Denote by $m$
the size of $P_A^{p-1}$ and $P_B^{q-2}$. By assumption, in $P_A^{p-1}$ --
hence
also in $P_B^{q-2}\sim P_A^{p-1}$ -- $x$ occupies $c(r-1,s-1)$ and
$c(r-1,s)$
is either empty or contains a region-3 element $z>x$.
But we saw above that in $P_B^{q-1}$, $c(r-1,s)$
contained $e\le x$. It follows that insertion step $P_B^{q-2}\rightarrow
P_B^{q-1}$ consisted of $e$ being inserted into $c(r-1,s)$, and -- if
$c(r-1,s)$ were previously occupied -- of bumping from it some $z>x$.
If $e$ bumped some $z$, then in
$P_B^{q-1}\rightarrow P_B^q$, $z$ would have been inserted into either row
$r$ or column
$s+1$. But we saw earlier that $x$ was bumped to row $r$ in $P_B^{q-1}
\rightarrow P_B^q$,
which implies that $z$ must have bumped $x$ during this step. This leads to
a contradiction, since $z$ could not have bumped $x<z$. Thus when $e$
occupied
$c(r-1,s)$ during $P_B^{q-2}\rightarrow P_B^{q-1}$, it did not
bump any element; the cell was previously empty. By occupying an empty cell,
$e$ increased the size of the tableau from $m$ to $m+1$: $|P_B^{q-2}|=m$,
$|P_B^{q-1}|=m+1$. Hence $|P_B^q|\ge m+1$.

On the other hand, $P_A^{p-1}$ is of size $m$, and during
$P_A^{p-1}\rightarrow
P_A^p$, $x$ was
bumped from row $r-1$ by some smaller element, so the size of the tableau
did
not change. It follows that $|P_A^p|=m<|P_B^q|$, contradicting
$P_A^p\sim P_B^q$.

It follows that $c(r-1,s)$ in $P_A^{p+1}$ contains an element $b\le u_j$,
and
this implies that when $u_j$ enters column $s$ in $P_A^{p+1}\rightarrow
P_A^{p+2}$, it settles in $c(r,s)$, to
the right of $x$. This completes the proof of~2.2.1.\medskip

{\it Proof of~2.2.2:}
By~2.2.1, $u_j$ settles in $c(r,s)$ in $P_A^{p+2}$. We show
that this implies that~2.12.2 for $P_A^{p+2}\sim P_B^{q+1}$ is satisfied.
In the diagrams below the proof, the cells outside of region 2
are marked with a $\star$.

Recall that by ``Case 2.2'',
$u_j$ is located in $c(r,s-1)$ in $P_A^p$, so $c(r,s-1)$ is part of a
connected
component of region~2. Let $\tau=$ the number of $t_i$'s and $\mu=$ the
number of $u_j$'s in this connected component. By~2.12.2 of $P_A^p\sim
P_B^q$,
region 2 of $P_B^q$ contains a corresponding connected component (in the
same cells as that of $P_A^p$), with $\tau\ \ t_i$'s and $\mu\ \ u_j$'s.

Since $u_j$ is located in $c(r,s-1)$ in $P_A^p$, by strict row inequality,
$c(r,s)$ is either empty or in region 3. Similarly,
since $x=t_i$ occupies $c(r,s-1)$ in $P_A^p\rightarrow P_A^{p+1}$,
$c(r-1,s-1)$ must contain a region-1 element in $P_A^{p+1}$ (and in
$P_A^p$).
\overfullrule=0pt
$$ P_A^p\ :\ \ \vcenter{\tabskip=0pt\offinterlineskip\halign
{\lft{$#$}&&\hquad\lft{$#$}\hquad&\lft{$#$}\cr
&&&&&$s$ \cr
&&&&&\vdots \cr
&&\multispan5{\hrulefill}\cr
&&\vl& \star &\vl&&\vl \cr
&&\multispan5{\hrulefill}\cr
r &\cdots&\vl& u_j &\vl& \star &\vl \cr
&&\multispan5{\hrulefill}\cr}} $$\medskip
\noindent Thus, if $c(r-1,s)$ is in region 2 in $P_A^p$ and $P_B^q$, then in
both
tableaux it is part of a connected component distinct from that of
$c(r,s-1)$.
This other component consists of $\tau'$ $t_i$'s and $\mu'$ $u_j$'s.
(If $c(r-1,s)$ is not in region 2, then we let $\tau'=\mu'=0$.)

Since $c(r,s)$ is either empty or in region 3 in $P_A^p$
(and thus in $P_B^q$), it follows that the
same is true for $c(r,s+1)$ and $c(r+1,s)$. By the assumption at the
beginning
of case 2, during $P_B^q\rightarrow
P_B^{q+1}$, $x=t_i$ bumps some region-3 element $z$
from $c(r,s)$, thereby adding $c(r,s)$ to region 2. But $c(r,s)$ is adjacent
to
both $c(r,s-1)$ and $c(r-1,s)$, so when it joins region 2, it combines their
respective connected components into a single larger one. Since neither
$c(r,s+1)$ nor $c(r+1,s)$ is in region 2, it follows that in $P_B^{q+1}$,
$c(r,s)$ becomes part of a connected component
of region 2, consisting of $\tau+\tau'+1\ \ t_i$'s and $\mu+\mu'\ \ u_j$'s.
\overfullrule=0pt
$$ P_B^{q+1}\ :\ \vcenter{\tabskip=0pt\offinterlineskip\halign
{\lft{$#$}&&\hquad\lft{$#$}\hquad&\lft{$#$}\cr
&&&&&$s$ \cr
&&&&&\vdots \cr
&&\multispan7{\hrulefill}\cr
&&\vl& \star &\vl&&\vl&&\vl \cr
&&\multispan7{\hrulefill}\cr
r &\cdots&\vl&&\vl& x &\vl& \star &\vl \cr
&&\multispan7{\hrulefill}\cr
&&\vl&&\vl& \star &\vl&&\vl \cr
&&\multispan7{\hrulefill}\cr}} $$\medskip

Similarly, during $P_A^p\rightarrow P_A^{p+1}\rightarrow P_A^{p+2}$, $x=t_i$
bumps $u_j$ from $c(r,s-1)$, and $u_j$ bumps $z$ from $c(r,s)$, so the only
change in the shape of region 2 in $P_A^{p+2}$
is the addition of $c(r,s)$. Thus $c(r,s)$
of $P_A^{p+2}$ is part of a connected component of region 2, also
containing $\tau+\tau'+1\ \ t_i$'s and $\mu+\mu'\ \ u_j$'s, and this
component
is identically shaped to the corresponding component of $P_B^{q+1}$, so~2.12
.2
of $P_A^{p+2}\sim P_B^{q+1}$ is satisfied. This completes the proof
of~2.2.2.
\medskip

{\it Proof of~2.2.3:}
Property~2.12.1 is satisfied for $P_A^{p+2}\sim P_B^{q+1}$, since
in both $P_B^q\rightarrow
P_B^{q+1}$ and $P_A^{p}\rightarrow P_A^{p+1}\rightarrow P_A^{p+2}$, region 1
is unchanged, and the only change in region 3 is the elimination of
$c(r,s)$.

Since the same element $z$ is bumped from $c(r,s)$ in
$P_A^{p+1}\rightarrow P_A^{p+2}$ and $P_B^q\rightarrow
P_B^{q+1}$,~2.12.3 is satisfied, and the proof of~2.2.3 is complete.\medskip

It follows that $P_A^{p+2}\sim P_B^{q+1}$.
\medskip

\item
Suppose that $x=t_a<t_i$.\medskip
Case 3.1: Row $r$ in $P_A^p$ terminates with $z\le x$. Thus
$z$ is in region 1, so by~2.12.1 of $P_A^p\sim P_B^q$, row $r$ in $P_B^q$
terminates with $z$. In such a case, $x$ is affixed to the end of row $r$
in both tableaux, so $P_A^{p+1}$ and $P_B^{q+1}$
have the same shape, and clearly satisfy~2.12.1 and~2.12.2.
Let $m$ denote the size of $P_A^{p+1}$ and $P_B^{q+1}$.
If $m=n$, which is the size of $P_A$ and $P_B$, then the insertion algorithm
terminates here. Otherwise, the next step is to begin $v_{m+1}$'s insertion
path by inserting $v_{m+1}$ into either
the first row or the first column in both tableaux. This verifies~2.12.3
and we have $P_A^{p+1}\sim P_B^{q+1}$.\medskip

Case 3.2: Row $r$ in $P_A^p$ (and hence in $P_B^q$) contains an element
greater than $x$. In
both tableaux, $x$ bumps from row $r$ the leftmost element greater than
itself. By~2.12.1 and~2.12.2 of $P_A^p\sim P_B^q$,
the same cell -- denoted $c(r,s)$ -- becomes occupied by $x$ in both
tableaux.
Thus if the element bumped by $x$ is identical in the two tableaux, then
$P_A^{p+1}\sim P_B^{q+1}$.

Suppose, however, that
$x$ bumps different elements from the cell(s) $c(r,s)$ of $P_A^p$ and
$P_B^q$.
By $P_A^p\sim P_B^q$, these must be $t_i$ and $u_j$.
Since $x=t_a$ occupies $c(r,s)$ in $P_A^{p+1}$, Lemma~2.3 implies that in
$P_A^{p-1}$, $x$
occupied $c(r-1,s')$ with $s\le s'$. Thus the $c(r-1,s)$ element $g$
in $P_A^{p-1}$ was $g\le x<u_j,t_i$, so $c(r-1,s)$ was a region-1 cell.
Since
$x$ was subsequently bumped from row $r-1$ by an element smaller than
itself,
it follows that $c(r-1,s)$ is a region-1 cell also in $P_A^p$ (and $P_B^q$).
Similarly, since $x=t_a<t_i$ settles in
$c(r,s)$ in $P_A^{p}\rightarrow P_A^{p+1}$,
$c(r,s-1)$ is a region-1 cell in $P_A^{p+1}$ (and $P_B^{q+1}$), and in
$P_A^p$
(and $P_B^q$).\medskip
\overfullrule=0pt
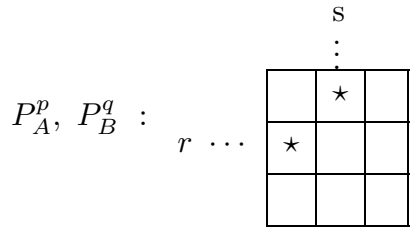
$$ P_A^p,\ P_B^q\ :\ \ \vcenter{\tabskip=0pt\offinterlineskip\halign
{\lft{$#$}&&\hquad\lft{$#$}\hquad&\lft{$#$}\cr
&&&&&$s$ \cr
&&&&&\vdots \cr
&&\multispan7{\hrulefill}\cr
&&\vl&&\vl& \star &\vl&\ &\vl \cr
&&\multispan7{\hrulefill}\cr
r &\cdots&\vl& \star &\vl&&\vl&\ &\vl \cr
&&\multispan7{\hrulefill}\cr
&&\vl&&\vl&&\vl&\ &\vl \cr
&&\multispan7{\hrulefill}\cr}} $$
$$\text{(The stars represent region-1 elements.)}$$
\medskip

Now, since $c(r,s)$ is in region 2 in both $P_A^p$ and $P_B^q$, by~2.12.2,
it
is part of a connected component of region 2 which is identically shaped,
and
contains the same number of $t_i$'s and $u_j$'s, in both tableaux. But
$c(r,s)$
contains $t_i$ in one tableau and $u_j$ in the other, so it follows that
at least one of $c(r,s+1)$ and $c(r+1,s)$ is in region 2 in $P_A^p$ and
$P_B^q$. By strict row and column inequality, this implies that $c(r,s)$
contains $t_i$ in $P_A^p$ and $u_j$ in $P_B^q$.

Denote by $C$ the connected component of region 2 containing $c(r,s)$.
Consider the subcomponent $C_1$,
consisting of all cells in $C$ which are to the right of or above $c(r,s)$.
In $P_A^p$, let $\alpha_A = \#t_i$'s, and $\beta_A = \#u_j$'s in $C_1$;
define
$\alpha_B$ and $\beta_B$ similarly in $P_B^q$. If $c(r,s+1)$ is not in
region
2, then $C_1$ is empty and $\alpha_A=\beta_A=\alpha_B=\beta_B=0$. On the
other
hand, if $C_1$ is nonempty, then by strict row and column inequality, every
northwest proper corner cell of $C_1$ contains $t_i$ in $P_A^p$ and $u_j$ in
$P_B^q$.\medskip
\overfullrule=0pt
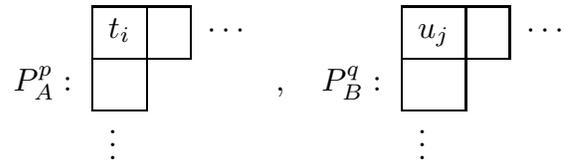
$$ P_A^p:\ \vcenter{\tabskip=0pt\offinterlineskip\halign
{\lft{$#$}&&\hquad\lft{$#$}\hquad&\lft{$#$}\cr
\multispan5{\hrulefill}\cr
\vl& t_i &\vl&\ &\vl&\cdots \cr
\multispan5{\hrulefill}\cr
\vl&&\vl \cr
\multispan3{\hrulefill}\cr
&\vdots \cr}}\ ,\quad
P_B^q:\ \vcenter{\tabskip=0pt\offinterlineskip\halign
{\lft{$#$}&&\hquad\lft{$#$}\hquad&\lft{$#$}\cr
\multispan5{\hrulefill}\cr
\vl& u_j &\vl&\ &\vl&\cdots \cr
\multispan5{\hrulefill}\cr
\vl&&\vl \cr
\multispan3{\hrulefill}\cr
&\vdots \cr}}\ $$\medskip

\noindent Similarly, every southeast proper corner cell of $C_1$ contains
$u_j$ in $P_A^p$ and $t_i$ in $P_B^q$.
\overfullrule=0pt
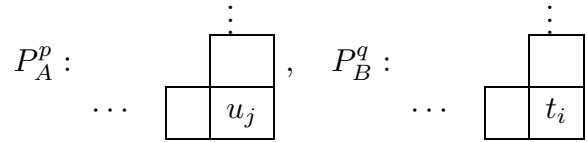
$$ P_A^p:\ \vcenter{\tabskip=0pt\offinterlineskip\halign
{\lft{$#$}&&\hquad\lft{$#$}\hquad&\lft{$#$}\cr
&&&&&\vdots \cr
&&&&\multispan3{\hrulefill}\cr
&&&&\vl&&\vl \cr
&&\multispan5{\hrulefill}\cr
\cdots&&\vl&\ &\vl& u_j &\vl \cr
&&\multispan5{\hrulefill}\cr}}\ ,\quad
P_B^q:\ \vcenter{\tabskip=0pt\offinterlineskip\halign
{\lft{$#$}&&\hquad\lft{$#$}\hquad&\lft{$#$}\cr
&&&&&\vdots \cr
&&&&\multispan3{\hrulefill}\cr
&&&&\vl&&\vl \cr
&&\multispan5{\hrulefill}\cr
\cdots&&\vl&\ &\vl& t_i &\vl \cr
&&\multispan5{\hrulefill}\cr}}\ $$\medskip

Consider the top row of $C_1$. If it contains more than one cell, then its
leftmost cell is a northwest corner. Thus the structure of $C_1$ is as in
the
following diagram, where for example, a cell marked $t_i/u_j$ contains $t_i$
in $P_A^p$ and $u_j$ in $P_B^q$. (A mark of $?$
denotes that a cell may contain either $t_i$ or $u_j$.)
and
elements
%
$$ {\aligned&\ \quad\quad\ t_i/u_j\cdots\cdots ?/t_i\\
&\quad\quad\quad\ \vdots\\ &\cdots\cdots u_j/t_i\endaligned}$$??

On the other hand, if
the top row of $C_1$ contains only one cell, then the structure of $C_1$ is
\medskip
$$ {\aligned&\quad\quad\quad\quad\quad\quad\quad\ \ ?/u_j\\
&\quad\quad\quad\quad\quad\quad\quad\quad\vdots\\
&\ \quad\quad\ t_i/u_j\cdots\cdots u_j/t_i\\
&\quad\quad\quad\ \vdots\\ &\cdots\cdots u_j/t_i\endaligned}$$




\noindent
In both cases it follows that $\alpha_B-\alpha_A=\beta_A-\beta_B\in\{0,1\}$.
We prove that

3.2.1: $\alpha_B-\alpha_A=\beta_A-\beta_B=1\implies P_A^{p+1}\sim
P_B^{q+2}$.

3.2.2: $\alpha_B-\alpha_A=\beta_A-\beta_B=0\implies P_A^{p+2}\sim
P_B^{q+1}$.
\medskip

Let $C_2$ be the subcomponent of $C$ consisting of all cells below or to
the left of $c(r,s)$. Let $\gamma_A=\# t_i$'s and $\delta_A=\# u_j$'s
in $C_2$ of $P_A^p$, and define $\gamma_B$ and $\delta_B$ similarly for
$P_B^q$. Since neither $c(r-1,s)$ nor $c(r,s-1)$ is in region 2, it follows
that $C=C_1+C_2+c(r,s)$. By~2.12.2 of $P_A^p\sim P_B^q$, $C$ contains the
same number of $t_i$'s and $u_j$'s in $P_A^p$ as in $P_B^q$.
In both tableaux, let $\tau$ be the number of $t_i$'s, and $\mu$ be
the number of $u_j$'s, in $C$. Since $c(r,s)$ contains $t_i$ in $P_A^p$ and
$u_j$ in $P_B^q$, it follows that:
$$ \tau=\alpha_A+\gamma_A+1=\alpha_B+\gamma_B,\quad
\mu=\beta_A+\delta_A=\beta_B+\delta_B+1. \leqno \star$$\medskip
{\it Proof of 3.2.1:}
Suppose that $\alpha_B-\alpha_A=\beta_A-\beta_B=1$.
Then $\gamma_A-\gamma_B=\delta_B-\delta_A=0$, so $C_2$ is either empty or
contains an equal number of $t_i$'s and $u_j$'s in $P_A^p$ as in $P_B^q$.
Also, $C_1$ is nonempty, so $c(r,s+1)$ is in region 2 in $P_A^p$ and in
$P_B^q$. In $P_B^q$, $c(r,s)$
contains $u_j$, so by strict row inequality $c(r,s+1)$ contains $t_i$. The
subsequent insertion steps are therefore\medskip
\def\toga{\hbox{\sevenrm $x$ bumps $t_i$\lower 12pt\hbox
{\hskip -40pt\sevenrm from $c(r,s)$}}}
\def\togb{\hbox{\sevenrm $t_i$ enters\lower 12pt\hbox{\hskip -30pt 
\sevenrm row $r+1$}}}
\def\togone{\hbox{\sevenrm $x$ bumps $u_j$\lower 12pt
\hbox{\hskip -45pt\sevenrm from $c(r,s)$}}}
\def\togtwo{\hbox{\sevenrm \hskip 5pt$u_j$ bumps $t_i$\lower 12pt
\hbox{\hskip -45pt\sevenrm from $c(r,s+1)$}}}
\def\togthr{\hbox{\sevenrm \hskip 10pt $t_i$ enters\lower 12pt
\hbox{\hskip -25pt\sevenrm row $r+1$}}}
$$ P_A^p~\ora_{\toga}~P_A^{p+1}~\ora_{\togb} $$\medskip
$$ P_B^q~\ora_{\togone}~P_B^{q+1}~\ora_{\togtwo}~P_B^{q+2}
\ora_{\togthr}~. $$
\medskip

Thus in both $P_A^p\rightarrow P_A^{p+1}$ and $P_B^q\rightarrow P_B^{q+1}
\rightarrow P_B^{q+2}$, $c(r,s)$ is eliminated from region 2, and we are
left
with two separate components $C_1$ and $C_2$ (and with $t_i$ to be inserted
into row $r+1$).
No change occurs in $C_2$, so in both $P_A^{p+1}$ and $P_B^{q+2}$, $C_2$ has
$\gamma_A=\gamma_B\ t_i$'s and $\delta_A=\delta_B\ u_j$'s. Similarly, no
change occurs in $C_1$ in $P_A^p\rightarrow P_A^{p+1}$, so $C_1$ of
$P_A^{p+1}$ contains $\alpha_A\ t_i$'s and $\beta_A\ u_j$'s. On the other
hand,
in $P_B^q\rightarrow P_B^{q+1}\rightarrow P_B^{q+2}$, a single change occurs
in $C_1$, when the $t_i$ in $c(r,s+1)$ is replaced with $u_j$. Thus $C_1$
of $P_B^{q+2}$ contains $\alpha_B-1\ t_i$'s and $\beta_B+1\ u_j$'s. But
by~3.2.1,
$\alpha_B-1=\alpha_A$ and $\beta_B+1=\beta_A$, so $C_1$ contains the same
number of $t_i$'s and $u_j$'s in $P_A^{p+1}$ as in $P_B^{q+2}$, and~2.12.2
is
satisfied for $P_A^{p+1}\sim P_B^{q+2}$.

Now in both $P_A^{p+1}$ and $P_B^{q+2}$,
the only change that occurs in region 1 is that the
same element $x$ is added to $c(r,s)$, so~2.12.1 is satisfied. Similarly, as
was already mentioned, both
$P_A^{p+1}\rightarrow P_A^{p+2}$ and $P_B^{q+2}\rightarrow P_B^{q+3}$
consist of $t_i$ entering row $r+1$, so~2.12.3 is satisfied. It follows
that $P_A^{p+1}\sim P_B^{q+2}$. This completes the proof of 3.2.1.
\medskip

{\it Proof of 3.2.2:}
The proof of 3.2.2 is dual, in a sense, to the proof of 3.2.1. Here are the
details.

Suppose that $\alpha_B-\alpha_A=\beta_A-\beta_B=0$.
Then $C_1$ is either empty or
contains an equal number of $t_i$'s and $u_j$'s in $P_A^p$ as in $P_B^q$.
Thus by ($\star$), $\gamma_B-\gamma_A=\delta_A-\delta_B=1$, so $C_2$ is
nonempty, which implies that $c(r+1,s)$ is in region 2 in $P_A^p$ and in
$P_B^q$. In $P_A^p$, $c(r,s)$
contains $t_i$, so by strict column inequality $c(r+1,s)$ contains $u_j$.
The
subsequent insertion steps are therefore\medskip
\def\togone{\hbox{\sevenrm $x$ bumps $t_i$\lower 12pt\hbox
{\hskip -40pt\sevenrm from $c(r,s)$}}}
\def\togtwo{\hbox{\sevenrm \hskip 10pt
$t_i$ bumps $u_j$\lower 12pt\hbox{\hskip -50pt 
\sevenrm from $c(r+1,s)$}}}
\def\togthr{\hbox{\sevenrm \hskip 18pt $u_j$ enters\lower 12pt
\hbox{\hskip -40pt\sevenrm column $s+1$}}}
\def\togfou{\hbox{\sevenrm \hskip -5pt
$x$ bumps $u_j$\lower 12pt
\hbox{\hskip -45pt\sevenrm from $c(r,s)$}}}
\def\togfiv{\hbox{\sevenrm \hskip 15pt $u_j$ enters\lower 12pt
\hbox{\hskip -40pt\sevenrm column $s+1$}}}
$$ P_A^p~\ora_{\togone}~P_A^{p+1}~\ora_{\togtwo}~P_A^{p+2}~
\ora_{\togthr} $$\medskip
$$ P_B^q~\ora_{\togfou}~P_B^{q+1}~\ora_{\togfiv} $$\medskip
Thus in both $P_A^p\rightarrow P_A^{p+1}\rightarrow P_A^{p+2}$ and
$P_B^q\rightarrow P_B^{q+1}$,
$c(r,s)$ is eliminated from region 2, and we are left
with two separate components $C_1$ and $C_2$ (and with $u_j$ to be inserted
into column $s+1$).
No change occurs in $C_1$, so in both $P_A^{p+2}$ and $P_B^{q+1}$, $C_1$ has
$\alpha_A=\alpha_B\ t_i$'s and $\beta_A=\beta_B\ u_j$'s. Similarly, no
change occurs in $C_2$ in $P_B^q\rightarrow P_B^{q+1}$, so $C_2$ of
$P_B^{q+1}$ contains $\gamma_B\ t_i$'s and $\delta_B\ u_j$'s. On the other
hand,
in $P_A^p\rightarrow P_A^{p+1}\rightarrow P_A^{p+2}$, a single change occurs
in $C_2$, when the $u_j$ in $c(r+1,s)$ is replaced with $t_i$. Thus $C_2$
of $P_A^{p+2}$ contains $\gamma_A+1\ t_i$'s and $\delta_A-1\ u_j$'s. But
3.2.2 and ($\star$) imply that
$\gamma_A+1=\gamma_B$ and $\delta_A-1=\delta_B$, so $C_2$ contains the same
number of $t_i$'s and $u_j$'s in $P_A^{p+2}$ as in $P_B^{q+1}$, and~2.12.2
is
satisfied for $P_A^{p+2}\sim P_B^{q+1}$.

Now in both $P_A^{p+2}$ and $P_B^{q+1}$,
the only change that occurs in region 1 is that the
same element $x$ is added to $c(r,s)$, so~2.12.1 is satisfied. Similarly, as
was already mentioned, both
$P_A^{p+2}\rightarrow P_A^{p+3}$ and $P_B^{q+1}\rightarrow P_B^{q+2}$
consist of $u_j$ entering column $s+1$, so~2.12.3 is satisfied. It follows
that $P_A^{p+2}\sim P_B^{q+1}$. This completes the proof of 3.2.2.
\bqed\endroster\enddemo\bigskip\bigskip

\noindent{\titl 3\ \ The Proof of Theorem 5}\bigskip

Here we prove, for example, Theorem~5(b). The proofs of parts (c) and (d)
of that theorem are similar.

Given $v\in a_{k,l}(n)$ and shuffle $A$, the (regular, dual)-$A$-RSK forms 
the tableau pair $(P^*,Q^*)=(P^*(v,A),Q^*(v,A))$ by applying the regular 
RSK to the $t_i$'s, and the
dual conjugate RSK to the $u_j$'s of $v$ under shuffle $A$. For simplicity,
we refer to this algorithm as the dual-$A$-RSK. As in
the $A$-RSK, $P^*$ is the insertion tableau, and $Q^*$ is the 
recording tableau of $v$ under $A$. Here $P^*$ is what we call a dual-$A$-SSYT;
that is, it is weakly $A$-increasing in rows, and strictly $A$-increasing in 
columns.\bigskip

\demo{\bf Example~3.1}
Let $k=2,\ l=1$, and $A: \ u_1 < u_2 < t_1 < t_2.$
Let 
$$ v=\Big{(}{{1\cdots\cdots 4}\atop{u_1,t_1,t_2,u_1}}\Big{)}~. $$
Then,
\overfullrule=0pt
$$ v~\ora_{\hbox{\sevenrm dual-A-RSK}}~
\vcenter{\tabskip=0pt\offinterlineskip\halign
{\lft{$#$}&&\hquad\lft{$#$}\hquad&\lft{$#$}\cr
\multispan3{\hrulefill}\cr
\vl& u_1 &\vl \cr
\multispan3{\hrulefill}\cr}}\quad\quad
\vcenter{\tabskip=0pt\offinterlineskip\halign
{\lft{$#$}&&\hquad\lft{$#$}\hquad&\lft{$#$}\cr
\multispan5{\hrulefill}\cr
\vl& u_1 &\vl& t_1 &\vl \cr
\multispan5{\hrulefill}\cr}}\quad\quad
\vcenter{\tabskip=0pt\offinterlineskip\halign
{\lft{$#$}&&\hquad\lft{$#$}\hquad&\lft{$#$}\cr
\multispan7{\hrulefill}\cr
\vl&u_1&\vl&t_1&\vl&t_2&\vl \cr
\multispan7{\hrulefill}\cr}}\quad\quad
\vcenter{\tabskip=0pt\offinterlineskip\halign
{\lft{$#$}&&\hquad\lft{$#$}\hquad&\lft{$#$}\cr
\multispan7{\hrulefill}\cr
\vl& u_1 &\vl& u_1 &\vl& t_2 &\vl \cr
\multispan7{\hrulefill}\cr
\vl& t_1 &\vl \cr
\multispan3{\hrulefill}\cr}}\  =\ P^*,$$
\noindent and
\overfullrule=0pt
$$ Q^*\ =\  \vcenter{\tabskip=0pt\offinterlineskip\halign
{\lft{$#$}&&\hquad\lft{$#$}\hquad&\lft{$#$}\cr
\multispan7{\hrulefill}\cr
\vl& 1 &\vl& 2 &\vl& 3 &\vl \cr
\multispan7{\hrulefill}\cr
\vl& 4 &\vl \cr
\multispan3{\hrulefill}\cr}}\bqed $$
\enddemo\bigskip

\proclaim{\bf Lemma~3.2}
Let $v\in a_{k,l}(n)$, $A\in I$ and
$$ v~\ora_{\hbox{\sevenrm A-RSK}}~(P,Q),\quad\quad 
v~\ora_{\hbox{\sevenrm dual-A-RSK}}~(P^*,Q^*).$$

\noindent If $v$ is non-repeating in its $u$-elements, then 
$P=P^*$ and $Q=Q^*$.
\endproclaim\bigskip

\demo{\bf Proof}
The $A$-RSK and the dual-$A$-RSK differ in only one rule: When some $u_j$
enters a column under the $A$-RSK, it bumps the first element $w_m$ such
that $w_m>u_j$ (or if no such $w_m$ exists, it settles at the end of the
column). On the other hand, under the dual-$A$-RSK, $u_j$ bumps the first
element $w_r$ such that $w_r\ge u_j$ (or settles at the end of the column).
But $u_j$ may appear only once in $w$, which implies that $w_r>u_j$, so this
step is the same as that of the $A$-RSK. The proof now follows.
\bqed\enddemo\bigskip

\definition{\bf Notation}
$v\in a_{k,l}(n)$ is said to be of type $(\alpha_1,...,\alpha_k;
\beta_1,...,\beta_l)$ if it is a permutation of 
$t_1^{\alpha_1}\cdots t_k^{\alpha_k}u_1^{\beta_1}\cdots u_l^{\beta_l}$.
\enddefinition\bigskip

\proclaim{\bf Lemma~3.3}
Let $v\in a_{k,l}(n)$ be of type $(\alpha_1,...,\alpha_k;\beta_1,...,\beta_l)$,
and denote $\beta=\sum_{i=1}^l\beta_i$. Then there exists $w\in 
a_{k,\beta}(n)$ such that \roster
\item The $u$-elements of $w$ are non-repeating.
\item For every shuffle $A$, if $v~{\dsp\ora_{\hbox{\sevenrm 
dual-A-RSK}}}~(P^*_v,Q^*_v)$, 
then there exists a corresponding shuffle $A'$ of the elements of 
$w$ such that $w~{\dsp\ora_{\hbox{\sevenrm dual-A-RSK}}}~
(P^*_w,Q^*_w)$, 
where $P^*_w$ is identical to $P^*_v$ but with every $v_i$ changed to $w_i$ for 
all $i\le n$. Consequently, {\rm sh}$(P^*_v)=\ ${\rm sh}$(P^*_w)$.\endroster
\endproclaim\bigskip

\demo{\bf Proof}
To avoid confusion between the elements of $v$ and of $w$, we let $u'_1,...,
u'_l$ denote the $u$-elements of $v$.

Form the sequence $w$ from $v$ as follows. 
Replace the $u'_1$'s in $v$ with $u_1,...,u_{\beta_1}$, moving from right
to left.
Replace the $u'_2$'s with $u_{\beta_1+1},...,u_{\beta_1+\beta_2}$, moving
from right to left.
Continue in this way until $u'_l$, and including the $u'_l$'s. 

Clearly the $u$-elements of $w$ are non-repeating, satisfying~3.3.1. 

Given some shuffle $A$ of the elements of $v$, define the shuffle $A'$ of the
elements of $w$ as follows. For every $i\in\{1,...,k\}$, 
$$\aligned
&t_i<_{_A}u'_1\implies t_i<_{_{A'}}u_1<_{_{A'}}\cdots <_{_{A'}}u_{_{\beta_1}}
\\
&t_i<_{_A}u'_2\implies t_i<_{_{A'}}u_{_{\beta_1+1}}<_{_{A'}}\cdots 
<_{_{A'}}u_{_{\beta_1+\beta_2}}\\ 
&\quad\quad\quad\vdots\\ 
&t_i<_{_A}u'_l\implies t_i<_{_{A'}}u_{_{\beta_1+\cdots+\beta_{l-1}+1}}<_{_{A'}}
\cdots <_{_{A'}}u_{_{\beta}},\\ \\
&\quad\quad\quad\quad\quad\quad\quad\quad 
u_1'<_{_A}t_i\implies u_1<_{_{A'}}\cdots<_{_{A'}}
u_{_{\beta_1}}<_{_{A'}}t_i\\
&\quad\quad\quad\quad\quad\quad\quad\quad\quad\quad\quad\vdots\\
&\quad\quad\quad\quad\quad\quad\quad\quad  
u_l'<_{_A}t_i\implies u_{_{\beta_1+\cdots+\beta_{l-1}+1}}<_{_{A'}}\cdots
<_{_{A'}}u_{_\beta}<_{_{A'}}t_i.
\endaligned$$\medskip

We compare the $A$-RSK insertion of the $v$'s with the $A'$-RSK insertion of
the $w$'s. Note that the shuffle $A$ and its derived shuffle $A'$
are similar in that $v_i<_{_A}v_j\implies w_i<_{_{A'}}
w_j$, but they differ in one fundamental way: For $i<j$ such that $v_i,v_j,
w_i$ and $w_j$ are $u$-elements, $v_i=_{_A}v_j\implies w_i>_{_{A'}}w_j$.
Now, if $w_j$ reaches a cell inhabited by $w_i>_{_{A'}}w_j$, then it bumps
$w_j$ to the next column, just as $v_j$ would bump $v_i=_{_A}v_j$ to the next
column under the dual-$A$-RSK. On the other hand, if $w_i$ reaches a cell 
inhabited by $w_j<_{_{A'}}w_i$, it settles below $w_j$, whereas $v_i$ would
bump $v_j=_{_A}v_i$ to the next column. However, such a situation never occurs,
since $i<j$ and $v_i=_{_A}v_j$ implies that every column reached by $w_j$ is 
first reached by $w_i$. The proof of this is as follows.

Suppose that for some $x$, $w_i=u_{x+1}$ and $w_j=u_x$. Then  
every column reached by $w_j$ is first reached by $w_i$, by induction on the
columns of $P^*_w$. Trivially, $w_i$ reaches column $1$ before $w_j$. By the 
induction assumption, $w_i=u_{x+1}$ is in column $c'$, $c'\ge c$. If $c'>c$, 
then we are done. Assume $c'=c': w_i=u_{x+1}$ is already in column $c$, and
$w_j=u_x$ is inserted into column $c$. It bumps the first $w_d$ such that
$w_d\ge w_j=u_x$. Now $v_i=_{_A}v_j$ implies that there does not exist any 
$t_z$ 
such that $w_j<_{_{A'}}t_z<_{_{A'}}w_i$. Hence $w_d=u_{x+1}=w_i$ is bumped
to column $c+1$.

This clearly extends to the general case $i<j,$ $v_i=v_j$, $w_i=u_y$, 
$w_j=u_x$, for general $y>x$.

Hence the steps of the dual-$A'$-RSK on $w$ are identical 
to the steps of the dual-$A$-RSK on $v$, but with every $v_i$, $i\le n$, 
changed to $w_i$. This implies that~3.3.2 is satisfied for $w$.
\bqed\enddemo\bigskip

\demo{\bf Example~3.4} Let $v=t_2u_2u_1u_1t_1$, and $A=t_1<t_2<u_1<u_2$. The
sequence $w=t_2u'_3u'_2u'_1t_1$ clearly satisfies~3.3.1; we show that it
satisfies~3.3.2 for shuffle $A$, by letting $A'=t_1<t_2<u'_1<u'_2<u'_3$. Under
shuffles $A$ and $A'$, 
 
$$v\ \overrightarrow{_{\ A\text{-RSK}\ }}\ (P^*_v,Q^*_v)\quad \text{and}\quad
w\ \overrightarrow{_{\ \text{dual-}A'\text{-RSK}\ }}\ (P^*_w,Q^*_w),$$
where
\overfullrule=0pt
$$ P^*_v\ = \  \vcenter{\tabskip=0pt\offinterlineskip\halign
{\lft{$#$}&&\hquad\lft{$#$}\hquad&\lft{$#$}\cr
\multispan9{\hrulefill}\cr
\vl& t_1 &\vl& u_1 &\vl& u_1 &\vl& u_2 &\vl \cr
\multispan9{\hrulefill}\cr
\vl& t_2 &\vl \cr
\multispan3{\hrulefill}\cr}}\ =\ 
\vcenter{\tabskip=0pt\offinterlineskip\halign
{\lft{$#$}&&\hquad\lft{$#$}\hquad&\lft{$#$}\cr
\multispan9{\hrulefill}\cr
\vl&v_5&\vl&v_4&\vl&v_3&\vl&v_2&\vl \cr
\multispan9{\hrulefill}\cr
\vl&v_1&\vl \cr
\multispan3{\hrulefill}\cr}}$$\medskip
\overfullrule=0pt
$$ P^*_w\ =\  \vcenter{\tabskip=0pt\offinterlineskip\halign
{\lft{$#$}&&\hquad\lft{$#$}\hquad&\lft{$#$}\cr
\multispan9{\hrulefill}\cr
\vl& t_1 &\vl& u'_1 &\vl& u'_2 &\vl& u'_3 &\vl \cr
\multispan9{\hrulefill}\cr
\vl& t_2 &\vl \cr
\multispan3{\hrulefill}\cr}}\ =\ 
\vcenter{\tabskip=0pt\offinterlineskip\halign
{\lft{$#$}&&\hquad\lft{$#$}\hquad&\lft{$#$}\cr
\multispan9{\hrulefill}\cr
\vl& w_5 &\vl& w_4 &\vl& w_3 &\vl& w_2 &\vl \cr
\multispan9{\hrulefill}\cr
\vl& w_1 &\vl \cr
\multispan3{\hrulefill}\cr}} $$\medskip
\noindent Thus~3.3.2 is satisfied for shuffle $A$.
\bqed\enddemo\bigskip

We can now give 
\bigskip

\demo{\bf The Proof of Theorem~5(b)}
Let $v$ be of type $(\alpha_1,...,\alpha_k;\beta_1,...,\beta_l)$, and
denote $\beta=\sum_{i=1}^l\beta_i$. 
Lemma~3.3 implies that there exists a sequence $w\in a_{k,\beta}(n)$ with
no repeating $u$-elements, and with
shuffles $A',B'$ such that

$$ w~\ora_{\hbox{\sevenrm dual-A'-RSK}}~(P^*_{A'},Q^*_{A'}),\quad\quad 
w~\ora_{\hbox{\sevenrm dual-B'-RSK}}~(P^*_{B'},Q^*_{B'}),$$ 
where sh$(P^*_{A'})=$ sh$(P^*_A)$ and sh$(P^*_{B'})=$ sh$(P^*_B)$. Since
$w$ contains no repetitions in its $u$-elements, Lemma~3.3 implies that 
$$ w~\ora_{\hbox{\sevenrm A'-RSK}}~(P^*_{A'},Q^*_{A'}),\quad\quad 
w~\ora_{\hbox{\sevenrm B'-RSK}}~(P^*_{B'},Q^*_{B'}).$$ Thus by 
Theorem~2, sh$(P^*_{A'})=$ sh$(P^*_{B'})$, which implies our result.
\bqed\enddemo\bigskip 

The proofs of parts (c) and (d) of Theorem~5 are similar to that of 
Theorem~5(b), since Lemma~3.3 can be applied also to the
(dual, regular)-$A$-RSK and the (dual, dual)-$A$-RSK. Both algorithms are 
$t$-dual; for simplicity, let $t_1',...,t_k'$ denote 
the $t$-elements of $v$. The $t$'s of the sequence $w$ of Lemma~3.3 for
parts (c) and (d) 
are set as follows. Replace the $t_1'$'s in $v$ with $t_1,...,t_{\alpha_1}$, 
moving from left to right. Replace the $t_2'$'s with $t_{\alpha_1+1},...,
t_{\alpha_1+\alpha_2}$, moving from left to right. Continue in this way until
$t_k'$, and including $t_k'$. 

Since the (dual, regular)-$A$-RSK of part (c) is $u$-regular, the $u$'s of $w$ 
are identical to those of $v$. However, the (dual, dual)-$A$-RSK of part (d)
is $u$-dual, so in this case the $u$'s of $w$ are derived the same way as in 
the proof of Lemma~3.3. Finally, shuffle $A'$ is derived from $A$ in parts (c) 
and (d) by methods analogous to that of part (b).\bigskip\bigskip

\noindent{\ti References}\bigskip

\noindent\ref\no 1\by A. Berele and A. Regev\paper Hook Young Diagrams with
Applications to Combinatorics and to Representations of Lie Superalgebras
\jour Adv. in Math.\vol 64 (2)\yr 1987\pages 118--175\endref
\ref\no 2\by A. Berele and J.B. Remmel\paper Hook Flag Characters and
Combinatorics\jour Jour. of Pure and Appl. Alg.\vol 35\yr 1985\pages
245\endref
\ref\no 3\by G. Olshanski, A. Regev and A. Vershik\paper Frobenius-Schur
Functions,~~{\it Studies in Memory of I. Schur} (Birkhauser),
to appear.\endref
\ref\no 4\by J. B. Remmel\paper The Combinatorics of $(k,l)$-hook Schur
Functions
\jour Comb. and Alg.\yr 1984 \pages 253-287\inbook Contemp. Math.\vol
34\publ AMS\endref
\ref\no 5\by J. B. Remmel\paper Permutation Statistics and $(k,l)$-hook
Schur Functions
\jour Discrete Math.\vol 67 (3)\yr 1987\pages 271-298\endref
\ref\no 6\by J. B. Remmel\paper A Bijective Proof of a Factorization Theorem
for $(k,l)$-hook Schur Functions\jour Linear and Multilinear Alg.\vol  28
(3)\yr 1990\pages 119-154\endref
\ref\no 7\by R. Stanley\inbook Enumerative Combinatorics\vol 2\publ
Cambridge Univ. Press\yr 1999\endref

\end